\RequirePackage{fix-cm}
\documentclass[smallextended]{svjour3}       % onecolumn (second format)
\smartqed  % flush right qed marks, e.g. at end of proof
\usepackage{graphicx}
\usepackage[latin1]{inputenc}
\usepackage{amsfonts}
\usepackage{amsmath}

\usepackage{color}

\begin{document}

\title{Distortion Representations of Multivariate Distributions%\thanks{Grants or other notes
}


\author{Jorge Navarro        \and
Camilla Cal\` i \and Maria Longobardi
\and Fabrizio Durante
}


\institute{J. Navarro \at
 Departamento Estadística e Investigación Operativa,
 Universidad de Murcia\\
 Campus de Espinardo, 30100 Murcia, Spain \\
%Tel.: +34 868883509\quad
 %            Fax: +34 868884182\\
 \email{jorgenav@um.es}\quad
ORCID: 0000-0003-2822-915X          
%  \\
           \and
C. Cal\` i, M. Longobardi (corresponding author) \at
Dipartimento di Biologia. Università di Napoli Federico II\\Via Cintia, I-80126 Napoli, Italy \\
\email{camilla.cali@unina.it} \quad ORCID: 0000-0002-7228-2399\\ 
\email{malongob@unina.it} \quad ORCID: 0000-0002-6772-5167
\and
F. Durante \at
Dipartimento di Scienze dell'Economia, Universit\`a del Salento\\ Lecce, Italy\\
\email{fabrizio.durante@unisalento.it}
\quad ORCID: 0000-0002-4899-1080
}

\date{Received: date / Accepted: date}

\maketitle

\begin{abstract}
The univariate distorted distribution were introduced in risk theory to represent changes (distortions) in the expected distributions of some risks. Later they were also applied to represent distributions of order statistics, coherent systems, proportional hazard rate (PHR) and proportional reversed hazard rate (PRHR) models, etc. In this paper we extend this concept to the multivariate setup. We show that, in some cases,  they are a valid alternative to the  copula representations especially 
when the marginal distributions may not be easily handled. Several relevant examples illustrate the applications of such representations in statistical modeling. They include the study of paired (dependent) ordered data, joint residual lifetimes, order statistics and coherent systems.
\end{abstract}

\section{Introduction}
%----------------------------------------------------------------------------
Distorted distributions were introduced in the theory of choice under risk (see \cite{LSS12,SCP16,W96,Y87}) to model the changes in the distribution of the risk variable under study. The distorted distribution is defined as $F_d(t)=d(F(t))$, where $F$ is the original distribution function and $d:[0,1]\to[0,1]$ is a distortion function (it is increasing, continuous and satisfies $d(0)=0$ and $d(1)=1$). However they can be applied in several contexts. For example, in reliability theory and survival analysis, they  can be used to represent the distributions of coherent systems with identically  distributed (ID) components  (see, e.g., \cite{NASS13}). This includes both the cases of independent identically distributed (IID) and dependent identically distributed (DID)  component lifetimes. In particular, they are also useful  to represent the distributions of  order statistics (i.e., the ordered data obtained from a sample) since they have the same distributions as $k$-out-of-$n$ systems. They were also used to define classes of prior distributions in Bayesian statistics (see \cite{ARS16}) and to represent conditional distributions (see \cite{ND17,NLP17,NS18}).

The distorted distributions were extended in \cite{NASS16} to represent univariate distributions as  distortions of $n\geq 2$ distribution functions. These representations were applied to study  the distribution of a single coherent system formed from $n$ components with different distributions.  They can also be used to represent ordered data from different populations (or in presence of outliers) and  to perform stochastic comparisons (see \cite{N18,NASS16,SS11}).

Several multivariate distortions have been proposed as well with the purpose of changing (shift) the distribution function of a given random vector $(X_1,\dots,X_n)$. For example, the distortion of the first kind proposed in Valdez and Xiao \cite{VX11} maintains the copula and distorts the marginals (see Section 2.1). Alternatively, the distortion of the third kind proposed there  maintains 
the marginals and replaces the copula by a distorted copula (see also \cite{DFS10,DS16,M05}). 
%However, these representations always lead to copula representations (based on the marginal distributions). 

Other authors propose alternative representations  for a given multivariate distribution $\mathbf{F}$ to the classical ones based on copulas  (see \cite{DS16,N06}). For example, Klüppelberg and Resnick \cite{KR08} proposed to use the Pareto-copula $C_P$ to represent $\mathbf{F}$. $C_P$ is a multivariate distribution having a common marginal standard  Pareto distribution $F_P(t)=1-1/t$ for $t\geq 1$. Hence 
$$C_P(v_1,\dots,v_n)=C(F_P(v_1),\dots,F_P(v_n))$$
for $v_1,\dots,v_n\in\mathbb{R}$ and
\begin{equation}\label{KR}
\mathbf{F}(x_1,\dots,x_n)=C_P(F_P^{-1}(F_1(x_1)),\dots, F_P^{-1}(F_n(x_n)) ),
\end{equation}
where $C$ is the copula of $ \mathbf{F}$, $F_1,\dots,F_n$ are its marginals  and $F_P^{-1}(u)=1/(1-u)$ for $u\in(0,1)$. Similar definitions can be proposed for other relevant distributions (normal, exponential, etc.).

%\quad

In this paper we introduce the concept of {\it multivariate distorted distribution} (MDD) extending the univariate concept given above. They provide alternative representations for a given multivariate distribution that can be represented as {\it distortions} of univariate distributions. These representations are similar to the classical copula representations. Actually, the copula representations are included in this general model. The main difference is that the MDD representation may not be built from the univariate marginals of the considered model, but from any set of univariate distributions.
%The main difference is that  multivariate distorted distributions are not necessarily built with the marginals. This fact provides a wide flexibility and allows us to obtain different (useful) representations. 
On the contrary, even for continuous distributions, these representations are not unique, that is, we do not have a result similar to Sklar's theorem for copulas.
We must note that the purpose of these representations is not to change (distort) the original distribution $\mathbf{F}$ but to provide alternative representations for it. Also note that the MDD representations defined here are different to the representations of type  \eqref{KR} obtained from a general continuous distribution function $G$ (see Section 2.1).

We provide several relevant examples were these representations are useful although other examples can be obtained as well. In the first example we provide a representation for the residual lifetimes of the working components in a system at a given time $t>0$ extending the results obtained in \cite{LP19,NLP17}. In the second one, we study ordered paired data obtaining a representation for  the joint distribution of the smallest  and the largest data. This representation can be used to estimate the largest order statistic from the smallest order statistic. This procedure can be applied, for instance, to study diseases of paired organs (eyes, kidneys, lungs, etc.).  The representation can be extended to the general case of  ordered data (order statistics) from a sample of dependent or independent  identically distributed random variables. Finally, we show that they can also be applied in Reliability Theory to represent the joint distribution of two different coherent systems based on the same components. In particular, this representation can be used to compute the system reliability and the expected system residual lifetime at the time of the first component failure. 

The rest of the paper is organized as follows. In the following Section 2 we define the multivariate distorted distributions obtaining their main properties.  The relevant examples are placed in Section 3. Some illustrations of simulated  ordered paired data sets are given in Section 4. The conclusions and open tasks for future research projects are  in Section 5.

\section{Multivariate distorted distributions}

Throughout the paper we use the terms `increasing' and `decreasing' in a wide sense, that is, they mean `non-decreasing' and `non-increasing', respectively. For example, a function $D:[0,1]^n\to[0,1]$ is increasing if $D(u_1,\dots,u_n)\leq D(v_1,\dots, v_n)$ whenever $0\leq u_i\leq v_i\leq 1$ for all $i=1,\dots,n$.

\subsection{Definition}

Let $\mathbf{X}=(X_1,\dots,X_n)$ be a random vector over a probability space $(\Omega,\mathcal{S},\Pr)$. 
Then the joint distribution function $\mathbf{F}$ of $\mathbf{X}$ is
$\mathbf{F}(x_1,\dots,x_n)=\Pr(X_1\leq x_1, \dots,X_n\leq x_n).$ 
The  (marginal) distribution of $X_i$ is 
$F_i(x_i)=\Pr(X_i\leq x_i)=\mathbf{F}(+\infty,\dots,+\infty, x_i,+\infty,\dots,+\infty)$
for $i=1,\dots,n$. It is well known that the probability in the $n$-dimensional rectangle (box) determined by the points $(x_1,\dots,x_n)$ and $(y_1,\dots,y_n)$ can be computed from $\mathbf{F}$ as
$$\Pr(x_1< X_1\leq y_1,\dots, x_n< X_n\leq y_n)=\triangle_{(x_1,\dots,x_n)}^{(y_1,\dots,y_n)} \mathbf{F},$$
where $x_i\leq y_i$ for $i=1,\dots,n$,
$$\triangle_{(x_1,\dots,x_n)}^{(y_1,\dots,y_n)} \mathbf{F}:=\sum_{z_i=x_i \text{ or } y_i} (-1)^{\mathbf{1}(z_1,\dots,z_n)} \mathbf{F}(z_1,\dots,z_n),$$
$\mathbf{1}(z_1,\dots,z_n)=\sum_{i=1}^n 1(z_i=x_i)$ and $1(A)=1$ (resp. $0$) if $A$ is true (false).

From Sklar's theorem (see, e.g., p. 42 in  \cite{DS16}), we know that $\mathbf{F}$ can be written as
\begin{equation}\label{C1}
\mathbf{F}(x_1,\dots,x_n)=C(F_1(x_1), \dots, F_n(x_n))
\end{equation}
for all $x_1,\dots,x_n$, where $F_1,\dots,F_n$ are the marginal distributions and $C$ is a {\it copula} function. Moreover, if all these marginal distributions are continuous, then the copula function $C$ is unique. For the basic properties of copulas  we refer the reader to  \cite{DS16,N06} and references therein. Now we are going to define a different concept that can also be used to represent $\mathbf{F}$.

\begin{definition}
	A multivariate distribution function  $\mathbf{F}$ is said to be a {\it multivariate distorted distribution} (MDD)  of the univariate distribution functions $G_1,\dots,G_n$ if there exists a continuous  function $D:[0,1]^n\to [0,1]$ such that 
	\begin{equation}\label{MDD1}
	\mathbf{F}(x_1,\dots,x_n)=D(G_1(x_1),\dots,G_n(x_n)) \text{ 	for all }x_1,\dots,x_n.
	\end{equation}
	We write $\mathbf{F}\equiv MDD(G_1,\dots,G_n)$, when $\mathbf{F}$ is a MDD of $G_1,\dots,G_n$.
\end{definition}

For instance, any multivariate distribution function  $\mathbf{F}$ is a MDD of its own marginals via the copula (see \eqref{C1}), that is, $\mathbf{F}\equiv MDD(F_1,\dots,F_n)$. But, for some choices of $G_1,\dots,G_n$,  $\mathbf{F}$ need not to be a MDD of $G_1,\dots,G_n$. For instance, take $\mathbf{F}$ being a bivariate Gaussian (normal) distribution and $G_1$ and $G_2$ being Bernoulli distribution functions of parameter $p$.

Note that we assume that $\mathbf{F}$ is a multivariate distribution function and that $D$ is a continuous function. The main properties are given in the next subsections. Several (useful) examples are included in Section 3.

\begin{remark}
	An alternative method not included in this general model is that of the representations based on G-copulas for a given continuous univariate distribution function $G$ similar to the Pareto-copulas proposed in \eqref{KR}. They are defined as follows: if $G$ is a continuous univariate distribution function and $C$ is a copula, the G-copula $C_G$ is defined as $C_G(v_1,\dots,v_n)=C(G(v_1),\dots,G(v_n))$
	for all $v_1,\dots,v_n\in \mathbb{R}$.
	First of all, we must say that $C_G$ need not be a copula in the usual sense since its common marginal distribution is $G$ (it is only a copula  when $G$ is the standard uniform distribution). Clearly, $C_G$ can be used to obtain the following representation for $\mathbf{F}$,
	\begin{equation}\label{KR-G}
	\mathbf{F}(x_1,\dots,x_n)=C_G(F_G^{-1}(F_1(x_1)),\dots, F_G^{-1}(F_n(x_n)) ).
	\end{equation}
	Note that $F_G^{-1}\circ F_i$, $i=1,\dots,n$ need not be distribution functions and so \eqref{KR-G} cannot be considered as a MDD representation. If all the marginals are continuous, from Sklar's theorem, the G-copula $C_G$ in \eqref{KR-G} is unique. 
\end{remark}

\begin{remark}
	The distortion of the first kind proposed in \cite{VX11} can be represented as MDD. They propose to change the original multivariate distribution $\mathbf{F}$ to the distorted one
	$$\mathbf{F}_{d_1,\dots,d_n}(x_1,\dots,x_n):=C(d_1(F_1(x_1))),\dots,d_n(F_n(x_n)))),$$
	for given univariate distortion functions $d_1,\dots,d_n$. In order to compare $\mathbf{F}$ and $\mathbf{F}_{d_1,\dots,d_n}$ we could represent this later function as a MDD with
	$$\mathbf{F}_{d_1,\dots,d_n}(x_1,\dots,x_n)=D(F_1(x_1),\dots, F_n(x_n)),$$
	where
	$D(u_1,\dots,u_n):=C(d_1(u_1),\dots,d_n(u_n))$
	for $u_1,\dots,u_n\in[0,1]$. In particular, if $d_1=\dots=d_n=d$ say, then $D$ is a d-copula.
\end{remark}

\subsection{Main properties}

The main properties of the distortion function $D$ are stated in the following proposition. The analogous result to Sklar's theorem for MDD is stated in item $(iii)$.

\begin{proposition}\label{prop2.1}
	Let $(X_1,\dots,X_n)$ be a random vector with  distribution function $\mathbf{F}$ that can be written as in \eqref{MDD1}. Then:
	\begin{itemize}
		\item[(i)] $D(u_1,\dots,u_{i-1},0,u_{i+1},\dots,u_n)=0$ for all $u_1,\dots,u_{i-1},u_{i+1},\dots,u_n\in [0,1]$. In particular, we have $D(0,\dots,0)=0$.
		\item[(ii)] $D(1,\dots,1)=1$. 
		\item[(iii)] If $G_1,\dots, G_n$ are continuous, then $D$ is unique and it can be extended to be a multivariate distribution function with support contained in $[0,1]^n$. In particular, $D$ is increasing in each variable.
		\item[(iv)] The $i$th marginal distribution of $D$ is $D_i(u)=F_i(G_i^{-1}(u))$ (where $G_i$ is the quasi-inverse function of $G_i$) and, in particular, $D$ is not necessarily a copula.
		%\item[(v)]    
	\end{itemize}
\end{proposition} 
%\begin{proof}
{\bf Proof.}	The properties in $(i)$ and $(ii)$  can be obtained by taking limits in \eqref{MDD1} for $x_i\to-\infty$ and  $\min(x_1,\dots,x_n)\to +\infty$, respectively. 
	
	To prove $(iii)$ let us consider $(u_1,\dots,u_n)$  with $0\leq u_i\leq 1$ for $i=1,\dots,n$. As $G_1,\dots, G_n$ are continuous distribution functions, there exists $(x_1,\dots,x_n)\in\mathbb{R}^n$ such that $u_i=G_i(x_i)$  for $i=1,\dots,n$. Then, we consider the random variables $V_i=G_i(X_i)$ for $i=1,\dots,n$, and we get
	\begin{align*}
	D(u_1,\dots,u_n)
	&=D(G_1(x_1),\dots,G_n(x_n))\\
	&=\mathbf{F}(x_1,\dots,x_n)\\
	&=\Pr(X_1\leq x_1,\dots,X_n\leq x_n)\\
	&=\Pr(G_1(X_1)\leq G_1(x_1),\dots,G_n(X_n)\leq G_n(x_n))\\
	&=\Pr(V_1\leq u_1,\dots,V_n\leq u_n)
	\end{align*}
	and so $D$ is the multivariate distribution function of $V_1,\dots,V_n$.
	In particular, $D$ is increasing in each variable.
	
To prove $(iv)$ for $i=1$, we note that
$$F_1(x_1)=F(x_1,+\infty,\dots,+\infty)=D(G_1(x_1),1,\dots,1)=D_1(G_1(x_1))$$
and so $D_1(u)=F_1(G_1^{-1}(u))$ holds for $u\in(0,1)$. The proof for $i=2,\dots,n$ is similar. 

To prove now that $D$ is not always a copula	we consider $G_i:=F_i^{\alpha_i}$ where $F_i$ are the marginals for $i=1,\dots,n$ and  $D(u_1,\dots,u_n):=C(u^{1/\alpha_1}_1,\dots,u^{1/\alpha_n}_n)$  for some $\alpha_1,\dots,\alpha_n>0$. In this case \eqref{MDD1} holds and  we get 
	$$D(u,1,\dots,1)=C(u^{1/\alpha_1},1,\dots,1)= u^{1/\alpha_1}\neq u$$
	whenever $u\in(0,1)$ and $\alpha_1\neq 1$. Therefore, $D$ is not a copula when $\alpha_1\neq 1$. \qed
%\end{proof}

\begin{remark}
	In particular, if the distribution functions $G_1,\dots,G_n$ are continuous,  since the random vector $(V_1,\dots,V_n)$ used in the preceding proof of $(iii)$ is a componentwise increasing transformation of $(X_1,\dots,X_n)$, then $(V_1,\dots,V_n)$ and $(X_1,\dots,X_n)$ have the same copula, that is,  $C_D=C$ and the respective Kendal's tau and Spearman's rho coefficients satisfy
	$\tau(\mathbf{F})=\tau(D)$ and $\rho(\mathbf{F})=\rho(D).$ 
	%Therefore $D$ contains all the information about the dependence structure of $\mathbf{F}$. 
\end{remark}

The converse of property $(iii)$ in the preceding proposition can be stated as follows. This is a very relevant property in order to built new multivariate probability models.

\begin{proposition}
	If $D$ is a continuous multivariate distribution function with support contained in $[0,1]^n$, then the function $\mathbf{F}$ defined by \eqref{MDD1} is a multivariate distribution function for all univariate distribution functions $G_1,\dots,G_n$. 
\end{proposition}
%\begin{proof}
{\bf Proof. }	Clearly, if \eqref{MDD1} holds for some distribution functions $G_1,\dots,G_n$, then% we have 
	$$\lim_{x_i\to -\infty}\mathbf{F}(x_1,\dots,x_n)=\lim_{x_i\to -\infty}D(G_1(x_1),\dots,G_n(x_n))
	%=D(G_1(x_1),\dots,G_{i-1}(x_{i-1}),0,G_{i-1}(x_{i-1}),\dots,G_n(x_n))
	=0$$
	for $i=1,\dots,n$ and 
\begin{align*}
\lim_{\min(x_1,\dots,x_n)\to +\infty}\mathbf{F}(x_1,\dots,x_n)
&=\lim_{\min(x_1,\dots,x_n)\to +\infty}D(G_1(x_1),\dots,G_n(x_n))\\
&=D(1,\dots,1)=1
\end{align*}
	since $D$ is continuous and the support of $D$ is included in $[0,1]^n$. Moreover, $\mathbf{F}$ is right-continuous in each variable since $D$ is continuous and $G_1,\dots,G_n$ are right-continuous.
	
	Let us consider  now $(x_1,\dots,x_n)\in\mathbb{R}^n$ and $(y_1,\dots, y_n)\in\mathbb{R}^n$ such that  $x_i\leq y_i$ for $i=1,\dots,n$. Then we define $u_i=G_i(x_i)$ and $v_i=G_i(y_i)$ for $i=1,\dots,n$. As $G_i$ is a distribution function, we have $0\leq u_i\leq v_i\leq 1$ for $i=1,\dots,n$. Therefore 
	$$\triangle_{(x_1,\dots,x_n)}^{(y_1,\dots,y_n)} \mathbf{F}=\triangle_{(u_1,\dots,u_n)}^{(v_1,\dots,v_n)} D\geq 0$$
	since $D$ is a multivariate distribution function. Therefore, $\mathbf{F}$ is a distribution function. \qed
%\end{proof}

\quad

The functions $D$ satisfying the properties stated in Proposition \ref{prop2.1} can be called {\it multivariate distortion functions}. They are just continuous distribution functions with supports contained in $[0,1]^n$. The set $\mathcal{C}_n$ of copulas  of dimension $n$ is a subset of the set $\mathcal{D}_n$ of multivariate distortion functions  of dimension $n$.  From the preceding proposition, if $D\in \mathcal{D}_n$, then the right hand side of \eqref{MDD1} always determines a multivariate distribution function. So we can obtain different multivariate models from $D$ just by changing the univariate distributions $G_1,\dots,G_n$.

\begin{remark}
	Notice that the distortion (or aggregation) functions considered in \cite{NASS16} to get univariate distribution functions from $n$ univariate distribution functions do not necessarily belong to $\mathcal{D}_n$ (i.e., in many cases, they are not multivariate distortion functions). For example, the aggregation function $Q(u_1,\dots, u_n)=(u_1+\dots +u_n)/n$ is not a multivariate distortion function since $Q(0,u_2,\dots, u_n)=(u_2+\dots +u_n)/n>0$ for all $u_2,\dots,u_n\in(0,1]$. 
\end{remark}

In the next proposition we show that if  \eqref{MDD1} holds, then a similar representation holds for the joint survival (or reliability) function
$\mathbf{\bar F}(x_1,\dots,x_n)=\Pr(X_1>x_1,\dots,X_n>x_n).$

\begin{proposition}
	If a multivariate distribution function $\mathbf{F}$ can be written as in \eqref{MDD1} for $D\in \mathcal{D}_n$, then the joint survival function  can be written as
	\begin{equation}\label{MDD2}
	\mathbf{\bar F}(x_1,\dots,x_n)=\hat D(\bar G_1(x_1),\dots,\bar G_n(x_n))
	\end{equation}
	for all $x_1,\dots,x_n$, where $\bar G_i=1-G_i$ is the survival function associated to $G_i$ and $\hat D\in \mathcal{D}_n$ is a multivariate distortion function determined by $D$.
\end{proposition}	

The proof is immediate since the probability $\Pr(X_1>x_1, \dots,X_n>x_n)$ can be obtained as the probability in the $n$-dimensional rectangle determined by $(x_1,\dots,x_n)$ and $(+\infty,\dots,+\infty)$. The distortion function $\hat D$ determined by this formula can be called {\it dual (or survival) distortion function} as in the univariate case. 
Note that $\hat D$ is determined by $D$ (and vice versa). In the univariate case we have $\hat D(u)=1-D(1-u)$. However, when $n>1$, 
$$\hat D(u_1,\dots,u_n)\neq 1- D(1-u_1,\dots,1-u_n).$$
The formula to get $\hat D$ from $D$  is similar to the expression for the survival copula $\hat C$ in term of the distributional copula $C$.  For example, if $n=2$, we have:
\begin{align*}
\mathbf{\bar F}(x_1,x_2)&=\Pr(X_1>x_1,X_2>x_2)\\
&=\triangle_{(x_1,x_2)}^{(+\infty,+\infty)} \mathbf{F}\\
&= \mathbf{F}(+\infty,+\infty)-\mathbf{F}(x_1,+\infty)-\mathbf{F}(+\infty,x_2)+\mathbf{F}(x_1,x_2)\\
&= 1-D(G_1(x_1),1)-D(1,G_2(x_2))+D(G_1(x_1),G(x_2))\\
&= 1-D(1-\bar G_1(x_1),1)-D(1,1-\bar G_2(x_2))+D(1-\bar G_1(x_1),1-\bar G(x_2))\\
&=\hat D(\bar G_1(x_1),\bar G_2(x_2) ),
\end{align*}
where $\hat D(u_1,u_2)= 1-D(1-u_1,1)-D(1,1-u_2)+D(1-u_1,1-u_2)$ for all $u_1,u_2\in[0,1]$.
Note that if $D\in\mathcal{D}_2$ is the distribution function of the random vector $(U_1,U_2)$, then $\hat D$ is the distribution function of $(1- U_1,1- U_2)$ whose support is included in $[0,1]^2$ as well. So $\hat{D}\in\mathcal{D}_2$. Representations \eqref{MDD1} and \eqref{MDD2} are equivalent but, sometimes, it could be better to use \eqref{MDD2} instead of \eqref{MDD1} (or vice versa). Some examples are given in the next section.

In the last property of this subsection we prove that, under mild conditions, we can choose a common distribution function $G$ instead of $G_1,\dots,G_n$.

\begin{proposition}
	If \eqref{MDD1} holds for continuous distribution functions $G_1,\dots,G_n$ that have  supports included in $S=(\alpha,\beta)$ (where $\alpha$ can be $-\infty$ and $\beta$ can be $+\infty$) and $G$ is a continuous distribution function with support $S$ which is strictly increasing  in $S$, then there exists a continuous distortion $D_G:[0,1]^n\to [0,1]$ such that
	\begin{equation}\label{MDDG}
	\mathbf{F}(x_1,\dots,x_n)=D_G(G(x_1),\dots,G(x_n)) \text{ 	for all }x_1,\dots,x_n.
	\end{equation}
	%	for all $x_1,\dots,x_n$.
\end{proposition}
%\begin{proof}
{\bf Proof.}	If $x_i\in S$, then $G_i(x_i)=G_i(G^{-1}(G(x_i)))$, where $G^{-1}:(0,1)\to S$ is the inverse function of $G$ in $S$. Therefore, \eqref{MDDG} holds for all $x_1,\dots,x_n\in S$  and for the continuous function defined by 
	$$ D_G(u_1,\dots,u_n):=D(G_1(G^{-1}(u_1)),\dots,   G_n(G^{-1}(u_n)))$$
	for $u_1,\dots,u_n\in (0,1)$. This function can be extended to $[0,1]^n$ (by using its continuity) and with this extension \eqref{MDDG} holds for any $x_1,\dots,x_n$ (since $D$ is continuous). \qed
%\end{proof}

\quad

In particular, if \eqref{MDD1} holds for nonnegative random variables $X_1, \dots,X_n$ with continuous  distributions, then we can choose a common standard exponential distribution with $G(t)=1-\exp(-t)$ for $t\geq 0$ in \eqref{MDD1}.

\subsection{Marginal distributions}

A relevant property of the representation \eqref{MDD1} is that all the marginal random vectors have also MDD. In particular, all the random variables $X_1,\dots,X_n$ have univariate distorted distributions from $G_1,\dots,G_n$, respectively. The results can be stated as follows. To simplify the notation we just consider marginals of the type $(X_1,\dots,X_m)$ for $1\leq m\leq n$. The expressions for the other marginals can be obtained in a similar way.

\begin{proposition}
	If \eqref{MDD1} holds and $1\leq m\leq n$, then the joint distribution function $\mathbf{F}_{1,\dots,m}$ of $(X_1,\dots,X_m)$ can be written as 
	\begin{equation}\label{MDD3}
	\mathbf{F}_{1,\dots,m}(x_1,\dots,x_m)=D_{1,\dots,m}(G_1(x_1),\dots,G_m(x_m))
	\end{equation}
	for all $x_1,\dots,x_m$, where
	$D_{1,\dots,m}(u_1,\dots,u_m):=D(u_1,\dots,u_m,1,\dots,1)$
	for all $u_1,\dots,u_m\in[0,1]$ and $D_{1,\dots,m}\in\mathcal{D}_m$.
\end{proposition}
%\begin{proof}
{\bf Proof. }	The joint distribution function of the marginal $(X_1,\dots,X_m)$ can be written as $\mathbf{F}_{1,\dots,m}(x_1,\dots,x_m)=\mathbf{F}(x_1,\dots,x_m,+\infty,\dots,+\infty)$
	for all $x_1,\dots,x_m$.
	Then \eqref{MDD3} is obtained from \eqref{MDD1} taking into account that $G_i(+\infty)=1$ for any distribution function $G_i$ and $i=m+1,\dots,n$. Finally, \eqref{MDD3} implies  $D_{1,\dots,m}\in\mathcal{D}_m$. \qed
%\end{proof}

\quad

In particular, the distribution function of $X_i$ can be written as
\begin{equation}\label{di}
F_i(x_i)=D(1,\dots,1, G_i(x_i),1,\dots, 1)=D_i(G_i(x_i))
\end{equation} 
for all $x_i$, where $D_i(u):=D(1,\dots,1,u,1,\dots,1)$ and the value $u$ is placed at the $i$th position. Clearly, we have $G_i=F_i$ when  
$D_i(u)=u$ for all $u\in[0,1]$, that is, when the $i$th marginal of $D$ is a uniform distribution over the interval $(0,1)$.
In particular, we obtain the copula representation when all the univariate marginals of $D$ have a standard uniform distribution.

\subsection{Probability density function and conditional distributions}

Let us assume in this subsection that $\mathbf{F}$ is absolutely continuous with joint probability density function (PDF) $\mathbf{f}=_{a.e.} \partial_{1,\dots,n}\mathbf{F}$, where $\partial_i \mathbf{F}$ represents the partial derivative of $\mathbf{F}$ with respect to its $i$th variable, $\partial_{i,j} \mathbf{F}:=\partial_i \partial_j \mathbf{F}$, and so on. Then the joint PDF of a multivariate distorted distribution can be obtained as follows.

\begin{proposition}
	If \eqref{MDD1} holds for absolutely continuous distribution functions $G_1,\dots, G_n$ and a distortion function $D$ that admits continuous mixed derivatives of order $n$, then a joint PDF $\mathbf{f}$ of $(X_1,\dots,X_n)$ is 
	\begin{equation}\label{MDD4}
	\mathbf{f}(x_1,\dots,x_n)=_{a.e.}g_1(x_1) \dots g_n(x_n)\ \partial_{1,\dots,n}D (G_1(x_1),\dots,G_n(x_n)),
	\end{equation}
	where $g_i=_{a.e.} G'_i$ for $i=1,\dots,n$.% and $d=\partial_{1,\dots,n}D$.
\end{proposition}

The proof is immediate from \eqref{MDD1} and $\mathbf{f}=_{a.e.} \partial_{1,\dots,n}\mathbf{F}$. Note that if $D\in \mathcal{D}_n$ and it is absolutely continuous, then $d:=_{a.e}\partial_{1,\dots,n}D$ is a PDF of $D$.

We can also prove that all the conditional distributions also have (when they exist) multivariate distorted distributions. To simplify the notation we just consider the conditional distribution $(X_2|X_1=x_1)$. The result in such a case can be stated as follows.

\begin{proposition}
	If \eqref{MDD1} holds for two absolutely continuous distribution functions $G_1$ and $G_2$ and a distortion function $D$  that admits continuous mixed derivatives of order $2$, then the distribution function $F_{2|1}$ of $(X_2|X_1=x_1)$ can be written as 
	\begin{equation}\label{MDD5}
	F_{2|1}(x_2|x_1)=D_{2|1} (G_2(x_2)|G_1(x_1))
	\end{equation}
	whenever  $\lim_{v\to 0^+}\partial_1D(G_1(x_1),v)=0$, where  $D_{2|1}$ is a distortion function  given by
	$$D_{2|1}(v|G_1(x_1)):=\frac{\partial_{1}D(G_1(x_1),v)}{\partial_1D(G_1(x_1),1)}$$ 
	for $0<v<1$ and $x_1$ such that $\partial_1D(G_1(x_1),1)>0$. 
\end{proposition}
%\begin{proof}
{\bf Proof.}	The conditional  PDF of $(X_2|X_1=x_1)$ can be written as 
	$$f_{2|1}(x_2|x_1)=\frac{\mathbf{f}(x_1,x_2)}{f_1(x_1)}$$
	for all $x_1,x_2$ such that $f_1(x_1)>0$. Then by using \eqref{MDD4} and that $f_1(x_1)=g_1(x_1)D'_1(G_1(x_1))>0$, where $D_1(u):=D(u,1)$ and $D'_1(u)=\partial_1 D(u,1)$, we obtain
	$$f_{2|1}(x_2|x_1)%=\frac{g_1(x_1)g_2(x_2)}{g_1(x_1)d'_1(G_1(x_1))}\partial_{1,2}D(G_1(x_1),G_2(x_2))
	=_{a.e.}g_2(x_2)\frac{\partial_{1,2}D(G_1(x_1),G_2(x_2))}{\partial_1D(G_1(x_1),1)}.$$
	Thus, the conditional distribution function can be obtained as
	$$F_{2|1}(x_2|x_1)=\int_{-\infty}^{x_2} f_{2|1}(z|x_1) dz=  \int_{-\infty}^{x_2} g_2(z)\ \frac{\partial_{1,2}D(G_1(x_1),G_2(z))}{\partial_1D(G_1(x_1),1)}dz.$$
	Now, if we assume $\lim_{v\to 0^+}\partial_1D(G_1(x_1),v)=0$, then 
	$$F_{2|1}(x_2|x_1)=  \left[\frac{  \partial_{1}D(G_1(x_1),G_2(z))}
	{\partial_1D(G_1(x_1),1)}\right]_{z=-\infty}^{x_2}=\frac{\partial_{1}D(G_1(x_1),G_2(x_2))}{\partial_1D(G_1(x_1),1)}.$$
	Hence,  \eqref{MDD5} holds. \qed
%\end{proof}

\quad

A similar expression can be obtained from $\hat D$ for the conditional survival function by using \eqref{MDD2}.
Note that the  PDF of $(X_2|X_1=x_1)$ can be written as 
\begin{equation}\label{MDD6}
f_{2|1}(x_2|x_1)=_{a.e.}g_2(x_2) d_{2|1} (G_2(x_2)|G_1(x_1)),
\end{equation}
where  
$$d_{2|1}(v|G_1(x_1)):=D^\prime_{2|1}(v|G_1(x_1))=\frac{\partial_{1,2}D(G_1(x_1),v)}{\partial_1 D(G_1(x_1),1)}$$ 
for $0<v<1$ (zero elsewhere). Hence, from \eqref{MDD6}, the regression curve to predict $X_2$ from $X_1$,  $m_{2|1}(x_1):=E(X_2|X_1=x_1)$, can be obtained as
$$m_{2|1}(x_1)=\int_{-\infty}^{+\infty} x_2g_2(x_2) d_{2|1} (G_2(x_2)|G_1(x_1))dx_2.$$
If $X_2\geq 0$, an alternative expression  can be obtained from the conditional survival function. Another option to predict $X_2$ from $X_1$ is to use the conditional median regression curve
$\tilde m_{2|1}(x_1):=F^{-1}_{2|1}(0.5|x_1)$
(see \cite{K05} or \cite{N06}, p.\ 217). Note that this function can be computed from \eqref{MDD5} as $F^{-1}_{2|1}(v|x_1)=G_2^{-1}(D^{-1}_{2|1}(v|G_1(x_1)))$ for $0<v<1$ if we are able to compute the inverse functions of $G_2$ and   $D_{2|1}(v|G_1(x_1))$. Moreover, we can obtain confidence bands in a similar way (see \cite{K05}).  For example, the $75\%$ and  $90\%$ quantile-confidence bands for $(X_2|X_1=x_1)$ are determined by $(F^{-1}_{2|1}(0.25|x_1),F^{-1}_{2|1}(0.75|x_1))$ and $(F^{-1}_{2|1}(0.05|x_1),F^{-1}_{2|1}(0.95|x_1))$. Some examples are given in Section 4.

\subsection{Stochastic comparisons}

Let us assume now that two random vectors  $\mathbf{X}=(X_1,\dots,X_n)$  and $\mathbf{Y}=(Y_1,\dots,Y_n)$ have multivariate distorted distributions with respective distortion functions $D_X$ and $D_Y$ and with the same baseline distribution functions $G_1,\dots,G_n$. Hence we have the following immediate results for the lower orthant $\leq_{lo}$ and upper orthant $\leq_{uo}$ orders. For the definitions and main properties of these stochastic orders see, e.g.,  \cite{SS07}, pages 308--314.

\begin{proposition}
Let $\mathbf{X}$  and $\mathbf{Y}$ have MDD from $G_1,\dots,G_n$ with respective distortion functions $D_X$ and $D_Y$. 
	\begin{itemize}
		\item[(i)] If $D_X\geq D_Y$, then $\mathbf{X}\leq_{lo}\mathbf{Y}$.
		\item[(ii)]If $\hat D_X\leq \hat D_Y$, then $\mathbf{X}\leq_{uo}\mathbf{Y}$.
	\end{itemize} 
\end{proposition}

Analogously, if they have the same distortion, we get the following results.

\begin{proposition}
	If $\mathbf{X}$  and $\mathbf{Y}$ have MDD with the same distortion function $D$ from  $G_1,\dots,G_n$ and $H_1,\dots,H_n$, respectively, and   $G_i\geq H_i$ holds  for $i=1,\dots,n$, then $\mathbf{X}\leq_{lo}\mathbf{Y}$ and $\mathbf{X}\leq_{uo}\mathbf{Y}$. 
\end{proposition}

The proof is immediate from \eqref{MDD1} and \eqref{MDD2} by noting that if they have the same distortion function $D$, then they  have the same dual distortion function $\hat D$ as well. Note that we can combine both propositions to compare MDD with different distortions and different univariate distributions.

\subsection{Relationships with copula representations}

We have seen in Proposition \ref{prop2.1}, $(iv)$, that some MDD representations can be obtained from the copula representation. Conversely, under some mild conditions, we can theoretically obtain the copula representation from the MDD representation (see below). However, we must say that in practice this is not always possible (see the examples in the following section).  In these cases the MDD representation can be useful.

We have already seen  that if \eqref{MDD1} holds, then from \eqref{di} the $i$th  marginal distribution function can be obtained as $F_i(x_i)=D_i(G_i(x_i))$, where $D_i$ is the $i$th marginal of $D$ for $i=1,\dots,n$. Let us see now how to determine the copula function $C$ from $D$.

\begin{proposition}
	If \eqref{MDD1} holds and the marginal distribution function $F_i$   is continuous and strictly increasing in its support $S_i=(\alpha_i,\beta_i)$ for $i=1,\dots,n$ (where $\alpha_i$ and $\beta_i$ can be $-\infty$ and $+\infty$), then 
	\begin{equation}\label{C2}
	\mathbf{F}(x_1,\dots,x_n)=C(F_1(x_1),\dots,F_n(x_n))
	\end{equation}
	is the unique copula representation for $\mathbf{F}$,	where 
	$$C(u_1,\dots,u_n):= D(G_1(F^{-1}_1(u_1)),\dots,G_n(F^{-1}_n(u_n)))$$
	for $u_1,\dots,u_n\in(0,1)$ and $F^{-1}_1,\dots,F^{-1}_n$ are the inverse functions of $F_1,\dots,F_n$ in $S_1,\dots,S_n$.
\end{proposition}
%\begin{proof}
{\bf Proof.}	From \eqref{MDD1}, we have 
\begin{align*}
\mathbf{F}(x_1,\dots,x_n)
&=D(G_1(x_1),\dots,G_n(x_n))\\
&=D(G_1(F^{-1}_1(F_1(x_1))),\dots,G_n(F^{-1}_n(F_n(x_n))) )
\end{align*}
	for all $x_i\in S_i$. Hence \eqref{C2} holds. Moreover, note that $C$ can be extended to $[0,1]^n$ to get the unique  copula  of $\mathbf{F}$ obtained from Sklar's theorem. \qed
%\end{proof}

\quad

Note that if $D$ has continuous marginal distributions $D_1,\dots, D_N$, then $C$ is also the (unique) copula copula of $D$. However,  in practice, it is not always possible to determine $F^{-1}_1,\dots,F^{-1}_n$ in closed forms as we can see in the following examples. In these cases, it could be better to use \eqref{MDD1} instead of \eqref{C2}.

\section{Relevant examples}

The purpose of this section is to show some examples where the MDD representations can be useful. More examples of this type could be also found.

\subsection{Joint residual lifetimes}

In this section we assume that $X_1,\dots,X_n$ represent the lifetimes of $n$ components in an engineering or biological system. So we assume that they are nonnegative almost surely. In these contexts, the univariate residual lifetimes $(X_i-t|X_i>t)$, $i=1,\dots,n$, at time $t>0$ play a relevant role. Their survival functions are
$$\bar{F}_{i,t}(x):=\Pr(X_i-t>x|X_i>t)=\frac{\bar{F}_i(t+x)}{\bar{F}_i(t)}$$ 
for all $x\geq 0$, whenever $\bar{F}_i(t)>0$ for  $i=1,\dots,n$.
For example, the mean residual lifetime (MRL)  function $m_i(t)=E(X_i-t|X_i>t)$ is used to define a stochastic order (the MRL order) and two aging classes (the increasing/decreasing MRL classes, denoted as IMRL and DMRL, respectively).

Analogously, if $X_1,\dots,X_n$ are dependent, then we can consider the MRL random vector
$$\mathbf{X}_t=(X_1-t,\dots,X_n-t|X_1>t,\dots,X_n>t)$$
for $t\geq0$ such that $\mathbf{\bar F}(t,\dots,t)=\Pr(X_1>t,\dots,X_n>t)>0$. Note that it is natural to consider a common time $t$ for the components. Here we just consider that, at a time $t$, all the components are working. Some results for the residual lifetime of the system under this assumption were obtain in \cite{N18,ND17}. We will consider other options later.

In the following proposition we prove that  $\mathbf{X}_t$ admits a MDD representation  for all $t\geq 0$. Specifically, we will obtain  the MDD representation in terms of $\bar{F}_{1,t},\dots,\bar{F}_{n,t}$ for the survival function  of $\mathbf{X}_t$, defined as
$$\mathbf{\bar F}_t(x_1,\dots,x_n)=\Pr(X_1-t>x_1,\dots,X_n-t>x_n|X_1>t,\dots,X_n>t)$$
for $x_1,\dots,x_n,t\geq 0$. The analogous expression for the joint distribution function can be obtained in a similar way.

\begin{proposition}
	If $\mathbf{\bar F}(t,\dots,t)>0$ for a $t\geq 0$, then 
	\begin{equation}\label{MDDt}
	\mathbf{\bar F}_t(x_1,\dots,x_n)=\hat D_t(\bar{F}_{1,t}(x_1),\dots,\bar{F}_{n,t}(x_n))
	\end{equation}
	for all $x_1,\dots,x_n\geq 0$ and  the following dual distortion function 
\begin{equation}\label{Dt}
\hat D_t(u_1,\dots,u_n):=
\frac{\hat C(\bar F_1(t) u_1,\dots,\bar F_n(t)u_n)}
{\hat C(\bar F_1(t),\dots,\bar F_n(t))},\  u_1,\dots,u_n\in[0,1],
\end{equation}
	which depends on $\bar{F}_1(t),\dots,\bar{F}_n(t)$.
\end{proposition}
%\begin{proof}
{\bf Proof. }	First we note that  $\Pr(X_i>t)\geq \Pr(X_1>t,\dots,X_n>t)>0$ for $i=1,\dots,n$. So we can consider the survival functions $\bar{F}_{1,t},\dots,\bar{F}_{n,t}$ of the marginal residual lifetimes at time $t$.
	Then we note that $\mathbf{\bar F}_t$ can be written as
	\begin{align*}
	\mathbf{\bar F}_t(x_1,\dots,x_n)
	&=\Pr(X_1-t>x_1,\dots,X_n-t>x_n|X_1>t,\dots,X_n>t)\\
	&=\frac{\Pr(X_1>t+x_1,\dots,X_n>t+x_n)}{\Pr(X_1>t,\dots,X_n>t)}\\
	&=\frac{	\mathbf{\bar F}(t+x_1,\dots,t+x_n)}{	\mathbf{\bar F}(t,\dots,t)}
	\end{align*}
	for $x_1,\dots,x_n\geq 0$. Now we use the following copula representation for 
	$\mathbf{\bar F}$ (obtained from Sklar's theorem)
	$	\mathbf{\bar F}(x_1,\dots,x_n)=\hat C(\bar F_1(x_1),\dots,\bar F_n(x_n)),$ 
	where $\hat C$ is a continuous survival copula of $\mathbf{\bar F}$. Hence
	\begin{align*}
	\mathbf{\bar F}_t(x_1,\dots,x_n)
	&=\frac{	\mathbf{\bar F}(t+x_1,\dots,t+x_n)}{	\mathbf{\bar F}(t,\dots,t)}\\
	&=\frac{\hat C(\bar F_1(t+x_1),\dots,\bar F_n(t+x_n))}
	{\hat C(\bar F_1(t),\dots,\bar F_n(t))}\\
	&=\frac{\hat C(\bar F_1(t)\bar F_{1,t}(x_1),\dots,\bar F_n(t)\bar F_{n,t}(x_n))}
	{\hat C(\bar F_1(t),\dots,\bar F_n(t))}\\
	&=\hat D_t(\bar{F}_{1,t}(x_1),\dots,\bar{F}_{n,t}(x_n))
	\end{align*}
and  \eqref{MDDt} holds for the  function $D_t$  \eqref{Dt}. Hence $\hat D_t\in\mathcal{D}_n$. \qed
%\end{proof}

\begin{remark}
Note that the survival functions $\bar{F}_{1,t},\dots,\bar{F}_{n,t}$ of the marginal residual lifetimes at time $t$ are not the marginal survival functions of the random vector $\mathbf{X}_t$. The $i$th marginal survival function $\bar H_{i,t}$ of $\mathbf{X}_t$ is
\begin{align*}
\bar H_{i,t}(x)
&=\Pr(X_i-t>x|X_1>t,\dots,X_n>t)
\\&=\frac{\Pr(X_1>t,\dots,X_{i-1}>t,X_i>t+x,X_{i+1}>t,\dots,X_n>t)}{\Pr(X_1>t,\dots,X_n>t)}\\&
=\frac{	\mathbf{\bar F}(t,\dots,t,t+x,t,\dots,t)}{	\mathbf{\bar F}(t,\dots,t)}\\&
=\frac{\hat C(\bar F_1(t),\dots,\bar F_{i-1}(t),\bar F_i(t+x),\bar F_{i+1}(t+x)\dots,\bar F_n(t))}{\hat C(\bar F_1(t),\dots,\bar F_n(t))}.
\end{align*}
Hence representation \eqref{MDDt} is not a copula representation. To obtain the survival copula $\hat C_t$ of $\mathbf{X}_t$ from \eqref{MDDt} we need the inverse functions of 
$\bar H_{1,t},\dots,\bar H_{n,t}$.  Note that in many models it is not possible to get these inverse functions in closed forms. However, representation \eqref{MDDt} always holds and it is based on the univariate residual survival functions.
\end{remark}

We can obtain (in a similar way) MDD representations for other residual lifetimes. For example, if we know that at time $t$, the first $n-1$ components are alive but the $n$th component has failed, then the survival function of the random vector $$\mathbf{X}^{(n)}_t:=(X_1-t,\dots,X_n-t|X_1>t,\dots,X_{n-1}>t,X_n\leq t)$$
defined for $t>0$ such that $\Pr(X_1>t,\dots,X_{n-1}>t,X_n\leq t)>0$, can be written as
\begin{align*}
\mathbf{\bar F}^{(n)}_t(x_1&,\dots,x_{n-1})\\
&=\Pr(X_1-t>x_1,\dots,X_{n-1}-t>x_{n-1}|X_1>t,\dots,X_{n-1}>t,X_n\leq t)\\
&=\frac{\Pr(X_1>t+x_1,\dots,X_{n-1}>t+x_{n-1},X_n\leq x_n)}{\Pr(X_1>t,\dots,X_{n-1}>t,X_n\leq t)}\\
&=\frac{\Pr(X_1>t+x_1,\dots,X_{n-1}>t+x_{n-1})}{\Pr(X_1>t,\dots,X_{n-1}>t,X_n\leq t)}\\&\quad
-\frac{\Pr(X_1>t+x_1,\dots,X_{n-1}>t+x_{n-1},X_n>t)}{\Pr(X_1>t,\dots,X_{n-1}>t,X_n\leq t)}\\
&=\frac{	\mathbf{\bar F}(t+x_1,\dots,t+x_{n-1},0)-	\mathbf{\bar F}(t+x_1,\dots,t+x_{n-1},t)}{\Pr(X_1>t,\dots,X_{n-1}>t,X_n\leq t)}\\
&=\frac{	\mathbf{\bar F}(t+x_1,\dots,t+x_{n-1},0)-	\mathbf{\bar F}(t+x_1,\dots,t+x_{n-1},t)}{\mathbf{\bar F}(t,\dots,t,0)-	\mathbf{\bar F}(t,\dots ,t)}\\
&=\frac{	\hat C(\bar F_1(t+x_1),\dots,\bar F_{n-1}(t+x_{n-1}), 1)}{\hat C(\bar F_1(t),\dots,\bar F_{n-1}(t),1)-	\hat C(\bar F_1(t),\dots,\bar F_{n}(t))}
\\&\quad 
-\frac{	\hat C(\bar F_1(t+x_1),\dots,\bar F_{n-1}(t+x_{n-1}), \bar F_n(t))}{\hat C(\bar F_1(t),\dots,\bar F_{n-1}(t),1)-	\hat C(\bar F_1(t),\dots,\bar F_{n}(t))}
\end{align*}
for $x_1,\dots,x_{n-1}\geq 0$. Hence it can be written as 
$$
\mathbf{\bar F}^{(n)}_t(x_1,\dots,x_{n-1})=D^{(n)}_t(\bar{F}_{1,t}(x_1),\dots,\bar{F}_{n-1,t}(x_{n-1})),
$$
where $\bar{F}_{1,t},\dots,\bar{F}_{n-1,t}$ are the survival functions of the univariate residual lifetimes,  
$$
D^{(n)}_t(\mathbf{u})=\frac{	\hat C(\bar F_1(t)u_1,\dots,\bar F_{n-1}(t)u_{n-1}, 1)-	\hat C(\bar F_1(t)u_1,\dots,\bar F_{n-1}(t) u_{n-1}, \bar F_n(t))}{\hat C(\bar F_1(t),\dots,\bar F_{n-1}(t),1)-	\hat C(\bar F_1(t),\dots,\bar F_{n}(t))}
$$
for $\mathbf{u}=(u_1,\dots,u_{n})\in[0,1]^n$ and $
D^{(n)}_t\in\mathcal{D}_{n-1}$.

Similar representations can be obtained from \eqref{MDDt} for 
$(X_1-t,\dots,X_{n-1}-t|X_1>t,\dots,X_n>t)$
(a marginal of $\mathbf{X}_t$) and 
$(X_1-t,\dots,X_{n-1}-t|X_1>t,\dots,X_{n-1}>t)$
(the $n-1$ dimensional case). More interestingly, we can compare these two random vectors with  $\mathbf{X}^{(n)}_t$ just by comparing their distortion functions since their representations are based on the same univariate survival functions $\bar{F}_{1,t},\dots,\bar{F}_{n-1,t}$.  This is not the case if we use copula representations since these random vectors  have different marginal distributions when $X_1,\dots,X_{n}$ are dependent.  Let us see an example.

\begin{example}
	If $(X_1,X_2,X_3)$ have a common marginal survival function $\bar F$ and a Farlie-Gumbel-Morgenstern (FGM) survival copula
	$$	\hat C(u_1,u_2,u_3)=u_1u_2u_3[1+\theta(1-u_1)(1-u_2)(1-u_3)]%, \quad -1 \leq \theta \leq 1
	$$
	for $-1\leq \theta\leq 1$, let us consider the residual lifetimes
	$$
	\textbf{X}_t := (X_1-t, X_{2}-t | X_1>t, X_2>t), 
	$$ 
	$$
	\textbf{X}_t^{(3)} := (X_1-t, X_2-t | X_1>t, X_{2}>t, X_3\leq t) 
	$$
	and
	$$
	\textbf{X}_t^{*} := (X_1-t,X_{2}-t | X_1>t, X_{2}>t,X_3>t).
	$$
	If $t>0$ and $k:=\bar F(t)>0$, the respective dual distortion functions are
	$$\hat D_t(u_1,u_2)=\frac{\hat C(ku_1,ku_2,1)} {\hat C(k,k,1)}=u_1u_2,$$
	$$\hat D^{(3)}_t(u_1,u_2)=\frac{\hat C(ku_1,ku_2,1)-\hat C(ku_1,ku_2,k)} {\hat C(k,k,1)-\hat C(k,k,k)}=u_1u_2\frac{1-\theta k (1-ku_1)(1-ku_2)}{1-\theta k(1-k)^2},$$
	and
	$$\hat D^{\ast}_t(u_1,u_2)=\frac{\hat C(ku_1,ku_2,k)} {\hat C(k,k,k)}=u_1u_2\frac{1+\theta (1-ku_1)(1-ku_2)(1-k)}{1+\theta(1-k)^3}.$$
	A straightforward calculation shows that $\hat D^{\ast}_t\leq \hat D_t\leq \hat D^{(3)}_t$ when $\theta \leq 0$ and that $\hat D^{\ast}_t\geq \hat D_t\geq \hat D^{(3)}_t$  when $\theta \geq 0$. Hence $\textbf{X}_t^{*}\leq_{uo}\textbf{X}_t\leq_{uo}
	\textbf{X}_t^{(3)}$ for all $t>0$, all $\bar F$ and all $\theta\leq0$. The reverse orderings hold when $\theta\geq0$.
\end{example}

The same happens for other conditional residual lifetimes.  For example, an analogous representation can be obtained for 
$(X_1-t,\dots,X_{n-1}-t|X_1>t,\dots,X_{n-1}>t,X_n=t_1)$
when $0<t_1<t$ by using the techniques used in \cite{ND17,NS18}.
The  comparisons of these random vectors can be used to study the effect of the information available at time $t$ in the residual lifetimes of the working components at this time. We can  study inactivity times as well by using the procedures introduced in  \cite{NC19,NLP17}.

\subsection{Ordered paired data}

In this section we assume that we have a training sample $(X_1,Y_1), \dots, (X_m,Y_m)$ from a random vector $(X,Y)$ with  absolutely continuous joint distribution function $\mathbf{F}$.  For example, they may represent disease lifetimes for paired organs (breast, lung, eyes, etc.). However, in practice, for other individuals, we may just know $L=\min(X,Y)$ (dependent censored data) and we want to estimate $U=\max(X,Y)$. In particular, we can estimate the regression curve $m(t)=E(U|L=t)$ or its conditional survival function $S(x|t)=\Pr(U>x|L=t)$ (that can be used to compute the median regression curve and its confidence bands). Note that the target random variable is $U$ (not $Y$). We will consider other options later.

To this end we assume that $X$ and $Y$ have a common absolutely continuous distribution function $F$ and an  absolutely continuous copula $C$. Hence 
$\mathbf{F}(x,y)=C(F(x), F(y))$. 
In some cases, we may also assume that $C$ is permutation symmetric. In this case, $(X,Y)$ is exchangeable (EXC), that is, $(X,Y)=_{st}(Y,X)$, where $=_{st}$ denotes equality in law (distribution). Note that both $F$ and $C$ can be estimated from the training sample by using empirical or kernel type estimators (see e.g. the survey in \cite{S17}).

We want to obtain a MDD representation for the random vector $(L,U)$. Its joint distribution function $\mathbf{G}(x,y)=\Pr(L\leq x,U\leq y)$  can be computed as
$$\mathbf{G}(x,y)=\Pr(U\leq y)=\Pr(X\leq y,Y\leq y)=C(F(y),F(y))$$
when $y\leq x$ and as
$$\mathbf{G}(x,y)=\Pr(L\leq x,U\leq y)=\Pr( (\{X\leq x\}\cup \{Y\leq x\})\cap \{X\leq y\}\cap \{Y\leq y\})$$
when $x< y$. Hence, by using the inclusion-exclusion formula, we get
\begin{align*}
\mathbf{G}(x,y)
&=\Pr(X\leq x,Y\leq y)+\Pr(X\leq y,Y\leq x)-\Pr(X\leq x,Y\leq x)\\
&=C(F(x),F(y))+C(F(y),F(x))-C(F(x),F(x))
\end{align*}
for $x< y$. Therefore, 
\begin{equation}\label{MDDex1}
\mathbf{G}(x,y)=D(F(x),F(y))
\end{equation}
for the following distortion function
\begin{equation}\label{D}
D(u,v):=\left\{\begin{array}{crr}
C(v,v) & \text{for}& v\leq u;\\
C(u,v)+C(v,u)-C(u,u)& \text{for}& u< v.\\
\end{array}%
\right.
\end{equation}
Note that $D$ is continuous. Then the marginal distributions of $(L,U)$ can be written as
$$G_1(x):=\Pr(L\leq x)=D(F(x),1)=D_1(F(x))$$
and
$$G_2(y):=\Pr(U\leq y)=D(1,F(y))=D_2(F(y)),$$
where $D_1(u)=D(u,1)=2u-C(u,u)$
and $D_2(u)=D(1,u)=C(u,u)$ for all $u\in[0,1]$ (see \eqref{di}). Note that $D$ is not a copula and that we do not use the marginals $G_1$ and $G_2$ of $\mathbf{G}$. For example, if $X$ and $Y$ are independent, then  
$D_1(u)=D(u,1)=2u-u^2\neq u$ and $D_2(u)=D(1,u)=u^2\neq u$ for all $u\in(0,1)$. Note that $D_1$ and $D_2$ are univariate distortion functions (for all $C$).
%From the results given in Section 2, to get the copula representation of $\mathbf{G}$, we need the inverse function of $F$. However we prefer to use representation \eqref{MDDex1}.

From \eqref{MDD5} and \eqref{MDDex1}, the distribution function of $(U|L=x)$ can be obtained  as
\begin{equation}\label{G21}
G_{2|1}(y|x)=D_{2|1} (F(y)|F(x))
\end{equation}
for $y\geq x$, where 
$$D_{2|1}(v|F(x)):=\frac{\partial_{1}D(F(x),v)}{\partial_{1}D(F(x),1)},$$
$\partial_{1}D(u,v)=0$ for $v<u$ and 
$$\partial_{1}D(u,v)=\partial_{1}C(u,v)+\partial_{2}C(v,u)-\partial_{1}C(u,u)-\partial_{2}C(u,u)$$
for $v> u$. Hence $\lim_{v\to 0^+}\partial_{1}D(u,v)=0$ for all $0<u<1$. In particular, in the EXC case, we have 
$\partial_{1}D(u,v)=2\partial_{1}C(u,v)-2\partial_{1}C(u,u)$
and in the IID case  $\partial_{1}D(u,v)=2(v-u)$ for  $u\leq v\leq 1$. In this last case, we get $D_{2|1}(v|F(x))=(v-F(x))/\bar F(x)$ for $F(x)\leq v\leq 1$.

In the general case, its PDF is $g_{2|1}(y|x)=f(y)\ d_{2|1} (F(y)|F(x))$ for $y\geq x$ (zero elsewhere), where $f=F'$,
$$d_{2|1}(v|F(x)):=D^\prime_{2|1}(v|F(x))=\frac{\partial_{1,2}D(F(x),v)}{\partial_1 D(F(x),1)},$$
where  $d(u,v):=\partial_{1,2}D(u,v)=c(u,v)+c(v,u)$ for $v\geq u$ (zero elsewhere) is the PDF of $D$ and $c(u,v)=\partial_{1,2}C(u,v)$ for $0\leq u,v\leq 1$ (zero elsewhere) is the PDF of $C$. 
For example, in the IID case, we get $d(u,v)=\partial_{1,2}D(u,v)=2$ for $0\leq u\leq v\leq 1$ (i.e. a uniform distribution over this triangle).

Hence, the regression curve $m_{2|1}(x):=E(U|L=x)$, can be obtained as
$$m_{2|1}(x)=\int_{x}^{+\infty} yf(y) d_{2|1} (F(y)|F(x))dy=\int_{x}^{+\infty} yf(y)\frac{\partial_{1,2}D(F(x),F(y))}{\partial_1 D(F(x),1)} dy.$$ 
In particular, in the IID case, we get
$$m_{2|1}(x)=\int_{x}^{+\infty} zf(z)\frac{2}{2(1-F(x))} dz=E(X|X>x),$$
that is, $E(U-x|L=x)=E(X-x|X>x)$ which is the MRL of $X$ (or $Y$). Therefore,  the residual lifetime of $U$ from $x=L$ is equal to the residual lifetime of a component when we know that it is alive at time $x$ and $m_{2|1}$ uniquely determines $F$. In particular, if $F$ is an exponential distribution, then $m_{2|1}(x)=x+E(X)$.  

An alternative expression for the regression curve $m_{2|1}$ can be obtained from the survival function of $(U|L=x)$ that can be written as
$\bar G_{2|1}(y|x)=\hat D_{2|1} (\bar F(y)|\bar F(x))$
for $y\geq x$, where 
$\hat D_{2|1}(v|\bar F(x)):=1- D_{2|1}(1-v| F(x))$. As $(U-x|L=x)\geq 0$, then 
$$m_{2|1}(x)=x+E(U-x|L=x)=x+\int_x^{+\infty} \bar G_{2|1}(y|x) dy.$$
For example, in the IID case, we have 
$$\hat D_{2|1}(v|\bar F(x))=1- D_{2|1}(1-v| F(x))=1-\frac{(1-v)-F(x)}{1-F(x)}=\frac{v}{\bar F(x)}$$
for $0\leq v\leq \bar F(x)$. Hence, 
$\bar G_{2|1}(y|x)=\hat D_{2|1} (\bar F(y)|\bar F(x))=\bar F(y)/\bar F(x)$
for $y\geq x$, that is,  $(U|L=x)=_{st}(X|X>x)$ and so $m_{2|1}(x)=E(X|X>x)$ (as above).

\begin{figure}[ptb]
	\begin{center}
		\includegraphics[scale=0.45]{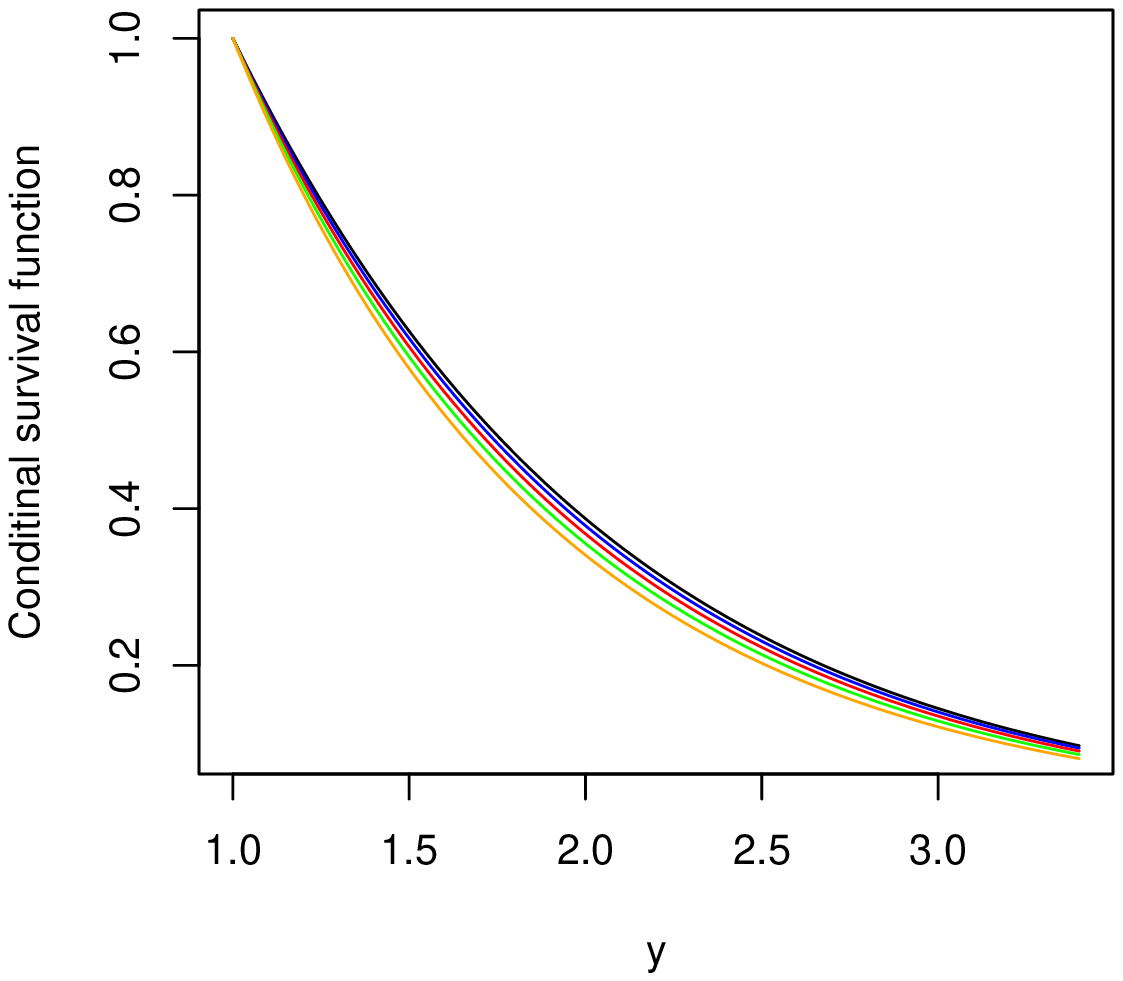}
		\includegraphics[scale=0.45]{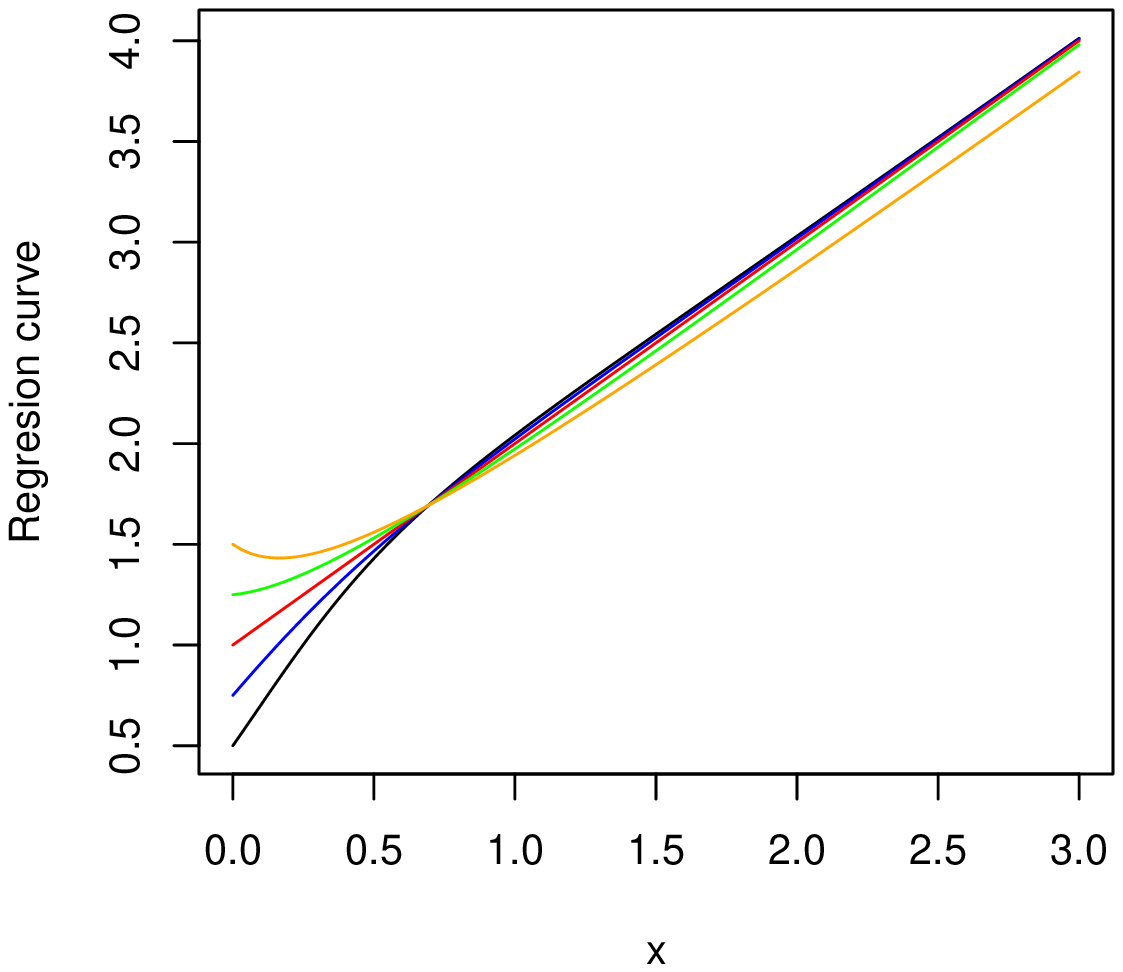}		%470x470despacho
	\end{center}
	\caption{Conditional survival functions $\bar G_{2|1}(y|1)$ (left) and  regression curves for $(L,U)$ when the marginals are  standard exponential distributions and $C$ is a FGM copula  with $\theta=-1,-0.5,0,0.5,1$ (orange, green, red, blue, black).} %The red curves represent the independent case ($\theta=0$).}%
	\label{figM}%
\end{figure}

This is true only for the IID case. For example,  let us consider the DID case with a common exponential distribution with mean $\mu=1$ and a FGM survival copula $\hat C(u,v)=uv[1+\theta(1-u)(1-v)]$ with $-1\leq \theta\leq 1$.	Then
$$\hat D_{2|1}(v|\bar F(x))=1- D_{2|1}(1-v| F(x))=\frac{v+\theta v(1-v)(1-2\bar F(x))}
{\bar F(x)+\theta \bar F(x)(1-\bar F(x))(1-2\bar F(x))}$$
for $0\leq v\leq \bar F(x)$ and 
\begin{align*}
m_{2|1}(x)&=x+\int_x^{+\infty}\bar G_{2|1}(y|x) dy\\
&=x+\int_x^{+\infty} \frac{ \bar F(y)+\theta  \bar F(y)(1- \bar F(y))(1-2\bar F(x))}
{\bar F(x)+\theta \bar F(x)(1-\bar F(x))(1-2\bar F(x))} dy
\\
&=x+ \frac{ (1+\theta-2\theta\bar F(x))\int_x^{+\infty} \bar F(y) dy -\theta  (1-2\bar F(x))\int_x^{+\infty} \bar F^2(y) dy}
{\bar F(x)+ \theta \bar F(x)(1-\bar F(x))(1-2\bar F(x))},
\end{align*}
where $\int_x^{+\infty} \bar F(y) dy=\bar F(x)E(X-x|X>x)$. If we assume now that $F$ is a standard exponential, then 
$$m_{2|1}(x)=x+\frac{1+\theta-2.5\theta e^{-x}+\theta e^{-2x}}
{1+\theta-3\theta e^{-x}+2\theta e^{-2x} }$$
for $x\geq 0$.  In Figure \ref{figM} we plot the conditional survival functions $\bar G_{2|1}(y|1)$ (left) and the regression curves for different values of $\theta$. In practice,  the dependence parameter $\theta$ of the FGM copula can be estimated from the Kendall's tau coefficient of the training sample (see \cite{N06}, p. 162) in order to estimate $U$ from $L=x$. A similar procedure can be applied to other families of copulas. If $F$ is an exponential distribution with mean $\mu>0$, then $E(U|L=x)=\mu E(U^*|L^*=x/\mu)$ for $x\geq 0$, where $L^*=L/\mu$ and $U^*=U/\mu$ are obtained from standard exponential distributions. Here  $\mu$ can be estimated from the training sample.

Another option to estimate $U$ from $L$ is to use the conditional median curve
$\tilde m_{2|1}(x):=G^{-1}_{2|1}(0.5|x)$
that can be obtained from \eqref{G21} and  the inverse function of  $D_{2|1}(u|F(x))$. Even more, we can obtain quantile-confidence bands (see Section 2.4). In the IID case, we get
\begin{equation}\label{IIDmedian}
\tilde m_{2|1}(x)=F^{-1}(F(x)+0.5\bar F(x)).
\end{equation}
In particular, if $F$ is an exponential distribution, then 
\begin{equation}\label{IIDmedianExp}\tilde m_{2|1}(x)=x-E(X)\ln(0.5)\approx x+0.6931472 E(X)<m_{2|1}(x)=x+E(X).
\end{equation}
An example with an EXC copula is given in Section 4.

Finally, the joint PDF of $(L,U)$ can be obtained from \eqref{MDD4} as
\begin{equation}\label{g2}
\mathbf{g}(x,y)= f(x)f(y)d(F(x),F(y))=f(x)f(y)\left[c(F(x),F(y))+c(F(y),F(x))\right]
\end{equation}
for $x\leq y$ (zero elsewhere). It can be used to get contour plots  (see Figure \ref{fig6}).

\subsection{Order statistics}

Let us consider now the ordered data $X_{1:n},\dots,X_{n:n}$ obtained from a sample $(X_1,\dots,X_n)$ of $n$ possibly dependent identically distributed (DID) random variables with absolutely continuous copula $C$ and marginal distribution $F$. The usual order statistics obtained from IID random variables are obtained when $C$ is the product copula. Recent results for conditional distributions in this case can be seen in \cite{AN20}. Clearly, the support of $(X_{1:n},\dots,X_{n:n})$ is included in the set $S=\{(x_1,\dots,x_n): x_1\leq \dots\leq x_n\}$. Then we can state the following result. 

\begin{proposition}
	The random vector $(X_{1:n},\dots,X_{n:n})$ have a MDD from a continuous distortion  function $D$ and  $F$, that is, its joint distribution function $\mathbf{G}$  can be written as
	$$\mathbf{G}(x_1,\dots,x_n)=D(F(x_1),\dots,F(x_n))\text{ for all }(x_1,\dots,x_n)\in S.$$
	%for all $x_1,\dots,x_n$. 
\end{proposition}

If  $F$  and $C$ are absolutely continuous with PDF $f$ and $c$,  then the PDF $\mathbf{g}$ of $\mathbf{G}$ is
$$\mathbf{g}(x_1,\dots,x_n)=f(x_1)\dots f(x_n)d(F(x_1),\dots,F(x_n))$$
for all $(x_1,\dots,x_n)\in S$ (zero elsewhere), where
$$d(u_1,\dots,u_n)=\sum_{\sigma\in \mathcal{P}_n}c(u_{\sigma(1)},\dots, u_{\sigma(n)})$$
for $0\leq u_1\leq \dots \leq u_n\leq 1$ (zero elsewhere),
$c=_{a.e.}\partial_{1,\dots,n}C$ is the PDF of $C$ and $\mathcal{P}_n$ is the set of permutations of dimension $n$.
If $C$ is EXC, then the expression of $d$ can be simplified to
$d(u_1,\dots,u_n)=n!\ c(u_{1},\dots, u_{n})$
for $0\leq u_1\leq \dots \leq u_n\leq 1$. 

The proof of the preceding proposition is straightforward but the explicit expression of $D$ is quite complicated.  If $n=2$ we obtain the distortion function of $(L,U)$ in the preceding subsection (see \eqref{D}). If $n=3$, the distortion function can be obtained as follows for $0\leq u_1\leq u_2 \leq u_3\leq 1$. If we assume $F(x)=x$ for $x\in[0,1]$, then 
\begin{align*}
D(u_1,u_2,u_3)&=\Pr(X_{1:3}\leq u_1,X_{2:3}\leq u_2,X_{3:3}\leq u_3)\\
&=\Pr( (A_1\cup A_2\cup A_3)\cap (A_{1,2}\cup A_{1,3}\cup A_{2,3})\cap A_{1,2,3}),
\end{align*}
where $A_i=\{ X_i\leq u_1\}$, $A_{i,j}=\{ X_i\leq u_2\}\cap \{ X_j \leq u_2\}$, and $A_{1,2,3}=\{ X_1\leq u_3\}\cap \{ X_2 \leq u_3\}\cap \{ X_3 \leq u_3\}$ for $i,j\in\{1,2,3\}$.  Hence 
\begin{align}\label{D3}
D(u_1,u_2,u_3)
&=\Pr( B_1\cup \dots  \cup B_9),
\end{align}
where $B_1=A_1\cap A_{1,2}\cap A_{1,2,3}$, $B_2=A_2\cap A_{1,2}\cap A_{1,2,3}$, $B_3=A_3\cap A_{1,2}\cap A_{1,2,3}$, 
$B_4=A_1\cap A_{1,3}\cap A_{1,2,3}$, $B_5=A_2\cap A_{1,3}\cap A_{1,2,3}$, $B_6=A_3\cap A_{1,3}\cap A_{1,2,3}$, $B_7=A_1\cap A_{2,3}\cap A_{1,2,3}$, $B_8=A_2\cap A_{2,3}\cap A_{1,2,3}$, and $B_9=A_3\cap A_{2,3}\cap A_{1,2,3}$.  Hence, the formula for $D$ is obtained by applying the inclusion-exclusion formula to \eqref{D3} taking into account that all these probabilities can be computed from $C$. For example
$$\Pr(B_1)=\Pr(A_1\cap A_{1,2}\cap A_{1,2,3})=\Pr(X_1\leq u_1,X_2\leq u_2,X_3\leq u_3)=C(u_1,u_2,u_3)$$
and
\begin{align*}
\Pr(B_1\cap B_2)
&=\Pr(A_1\cap A_2\cap  A_{1,2}\cap A_{1,2,3})\\%
&=\Pr(X_1\leq u_1,X_2\leq u_1,X_3\leq u_3)\\
&=C(u_1,u_1,u_3)
\end{align*}
$0\leq u_1\leq u_2 \leq u_3\leq 1$. The other probabilities can be obtained in a similar way. Clearly, this procedure can also be applied to the $n$ dimensional case  (but the expression for $D$ gets really involved).

\subsection{Coherent systems}

A {\it system} is a Boolean function $\psi:\{0,1\}^n\to\{0,1\}$ where $\psi(x_1,\dots,x_n)$ represents the state of the system that is completely determined by the components' states $x_1,\dots,x_n$.
Here $\psi(x_1,\dots,x_n)=1$ means that the system works and $\psi(x_1,\dots,x_n)=0$ that it has failed. A system is {\it semi-coherent} if $\psi$ is increasing, $\psi(0,\dots,0)=0$ and $\psi(1,\dots,1)=1$. It is {\it coherent} if $\psi$ is increasing and it is strictly increasing in every variable in at least a point. This last condition means that all the components are relevant for the system and it implies $\psi(0,\dots,0)=0$ and $\psi(1,\dots,1)=1$. For the basic properties of systems we refer the reader to the classic book \cite{BP75}. In particular, we can see that 
$\psi(x_1,\dots,x_n)=\min_{i=1,\dots,s}\max_{j\in \mathcal{C}_i}x_j,$ 
where $\mathcal{C}_1,\dots,\mathcal{C}_s$ are the minimal cut sets of the system. A set $A\subseteq\{1,\dots,n\}$ is a {\it cut set} of $\psi$ if $\psi(x_1,\dots,x_n)=0$ when $x_j=0$ for all $i\in A$.  A cut set is {\it minimal} if it does not contain other cut sets. The preceding expression can be used to extend $\psi$ to $\mathbb{R}^n$ and then the random system lifetime $T$ can be obtained as
$$T=\psi(X_1,\dots,X_n)=\min_{i=1,\dots,s}\max_{j\in \mathcal{C}_i}X_j,$$
where $X_1,\dots,X_n$ represent the random component lifetimes. Note that  $X^{\mathcal{C}_i}:=\max_{j\in \mathcal{C}_i}X_j$ represents the lifetime of the parallel system formed with the components in the set $\mathcal{C}_i$.

It is well known (see, e.g., \cite{N18,NASS16}) that the distribution function of a coherent (or semi-coherent) system $T$ can be obtained as a distortion of the components' distribution functions. The joint distribution of two semi-coherent systems with some shared components was studied in  \cite{NB10,NSB10}. Connections between dependence (copula) properties and ordering properties of systems were studied in \cite{NDF20}. 

In this section we consider two semi-coherent systems based on the same components with lifetimes $(X_1\dots, X_n)$ having a copula $C$ and a common marginal distribution function $F$. Then we can prove that the joint distribution of these two systems is a MDD. The result can be stated as follows. The proof is similar to that of Lemma 3.1 in \cite{NB10}, that is stated only for systems with IID components.

\begin{proposition}
	Let $T=\psi(X_1,\dots,X_n)$ and $T^*=\psi^*(X_1,\dots,X_n)$ be the lifetimes of two semi-coherent systems with  components' lifetimes $(X_1,\dots,X_n)$ having a copula $C$ and a common  marginal distribution function $F$. Then the joint distribution function $\mathbf{F}$ of $(T,T^*)$ can be written as
	\begin{equation}\label{T1T2}
	\mathbf{F}(x,y)=D(F(x),F(y))
	\end{equation}
	for all $x,y$ where $D\in \mathcal{D}_2$.
\end{proposition}

%\begin{proof}
{\bf Proof. }	Let $\mathcal{C}_1,\dots, \mathcal{C}_s$ and $\mathcal{C}^*_1,\dots, \mathcal{C}^*_{s^*}$ be the minimal cut sets of $T$ and $T^*$, respectively.  Hence the joint distribution $\mathbf{F}(x,y)=\Pr(T\leq x, T^*\leq y)$  can be written as
	\begin{align*}
	\mathbf{F}(x,y)
	&=\Pr(\min_{i=1,\dots,s}\max_{k\in \mathcal{C}_i}X_k \leq x, \min_{j=1,\dots,s^*}\max_{k\in \mathcal{C}^*_j}X_k \leq y)\\
	&=\Pr((A_1\cup\dots  \cup A_s)\cap (A^*_1\cup\dots\cup  A^*_{s^*}))\\
	&=\Pr\left( \cup_{i=1}^s\cup_{j=1}^{s^*} B_{i,j} \right),
	\end{align*}
	where $A_i:=\{\max_{k\in \mathcal{C}_i}X_k \leq x\}$, $A^*_j:=\{\max_{k\in \mathcal{C}^*_j}X_k \leq y\}$ and $B_{i,j}:=A_i\cap A^*_j$. Now we can apply the inclusion-exclusion formula to the union of the sets $B_{i,j}$.  Moreover we note that  if $x\leq y$, then
	\begin{align*}
	\Pr(B_{i,j})&
	=\Pr\left( \max_{k\in \mathcal{C}_i} X_k \leq x,  \max_{k\in \mathcal{C}^*_j}X_k \leq y \right)\\
	&=\Pr\left( \max_{k\in \mathcal{C}_i}X_k \leq x,  \max_{k\in \mathcal{C}^*_j-\mathcal{C}_i}X_k \leq y \right)\\
	&=C_{\mathcal{C}_i,\mathcal{C}^*_j}(F(x),F(y)),
	\end{align*}
	where $\mathcal{C}^*_j-\mathcal{C}_i
	=\mathcal{C}^*_j\cap  \mathcal{\bar C}_i$ 
	($\bar A$ is the complementary set of the set $A$), $C_{\mathcal{C}_i,\mathcal{C}^*_j}(u,v):=C(u_1,\dots,u_n)$, $u_k=F(x)$ if $k\in \mathcal{C}_i$, $u_k=F(y)$ if $k\in \mathcal{C}^*_j-\mathcal{C}_i$, and $u_k=1$ if $k\notin \mathcal{C}_i\cup \mathcal{C}^*_j$. Similar expressions can be obtained for the other probabilities in the inclusion-exclusion formula as $\Pr(B_{i,j}\cap B_{\ell,r})$, $\dots$ and for   $x>y$. Hence, we obtain  \eqref{T1T2}. \qed
%\end{proof}

\quad

Note that the distortion function $D(u,v)$ has different expressions for $u\leq v$ and for $u>v$. Moreover, note that $\mathbf{F}$ can have a singular part since the systems can fail at the same time (even if the joint distribution of the component lifetimes $(X_1,\dots,X_n)$ is absolutely continuous). 
Similar results can be stated for the joint distribution and the joint survival functions of $m$ coherent (or semi-coherent) systems based on the same components. 

In practice, the above result will be typically applied to study the joint distribution of a system $T^*$ and a related system $T<T^*$. The most usual situation is to consider $T=X_{1:n}=\min(X_1,\dots,X_n)$. In this case we want to study (predict) the system lifetime $T^*$ when we know the first component failure $T$. Note that we can use here all the results given in Section 2. 

For example, let us  consider the system $T^*=\max(X_1,\min(X_2,X_3))$ and $T=X_{1:3}$. If the joint distribution of $(X_1,X_2,X_3)$ is absolutely continuous, then $T^*<T$. Hence the joint survival function $\mathbf{\bar F}$ of $(T,T^*)$ can be written as 
$\mathbf{\bar F}(x,y)=\bar D(\bar F(x),\bar F(y))$
for all $x\leq y$, where  
$\bar D(u,v)=\hat C(u,v,v)+\hat C(v,u,u)-\hat C(v,v,v)$
for $0\leq v\leq u\leq 1$. 

Another interesting example is to consider the parallel  system with three components $T^*=X_{3:3}=\max(X_1,X_2,X_3)$ and its first components' failure (series system) $T=X_{1:3}$. Note that the distribution $\mathbf{F}$ of $(T,T^*)$ is a marginal of the random vector $(X_{1:3},X_{2:3},X_{3:3})$ (studied in the preceding subsection).  Note that $\mathbf{F}$ can be written as \eqref{T1T2} for  
\begin{align*}
D(u,v)&=C(u,v,v)+C(v,u,v)+C(v,v,u)-C(u,u,v)-C(u,v,u)\\
&\quad -C(v,u,u)+C(u,u,u)
\end{align*}
for $0\leq u\leq v\leq 1$ and $D(u,v)=C(v,v,v)$ for $0\leq v\leq u\leq 1$. If $C$ is EXC, then %the first expression can be reduced to
$$D(u,v)=3C(u,v,v)-3C(u,u,v)+C(u,u,u)$$
for $0\leq u\leq v\leq 1$. In particular, if the components are IID, we get
$$D(u,v)=3uv^2-3u^2v+u^3$$
for $0\leq u\leq v\leq 1$. Note that $T$ and $T^*$ are dependent even if $X_1,X_2,X_3$ are IID.

These expressions can be used to compute the joint PDF, the marginal and conditional distributions and the conditional median (jointly with the associated confidence regions). For example, in the last IID  case, the joint PDF is 
$$\mathbf{g}(x,y)=f(x)f(y)d(F(x),F(y))$$
where $d(x,y)=\partial_{1,2}D(u,v)=6(v-u)$ 
for $0\leq u\leq v\leq 1$ (zero elsewhere). The marginal distributions are 
$G_1(x)=\Pr(X_{1:3}\leq x)=D_1(F(x))$
and 
$G_2(y)=\Pr(X_{3:3}\leq y)=D_2(F(y)),$
where $D_1(u)=D(u,1)=3u-3u^2+u^3$ and $D_2(u)=D(1,u)=u^3$ for $u\in [0,1]$. Analogously, from \eqref{MDD5}, the conditional distribution function of $T$ given $T^*$ can be written as 
$$G_{2|1}(y|x)=\Pr(X_{3:3}\leq y|X_{1:3}=x)=D_{2|1} (F(y)|F(x)),$$
where  $D_{2|1}$ is given by
$$D_{2|1}(v|u):=\frac{\partial_{1}D(u,v)}{\partial_1D(u,1)}=\frac{3v^2-6uv+3u^2}{3-6u+3u^2}$$ 
for $0<u\leq v\leq 1$ since $\lim_{v\to 0^+}\partial_1D(u,v)=0$ for all $u\in(0,1)$.

\section{An illustration about paired ordered data}

We consider the problem stated in Section 3.2. We assume here that we have a sample $(X_1,Y_1), \dots, (X_m,Y_m)$ from a random vector $(X,Y)$ with joint absolutely continuous  distribution function.   We also assume that $X$ and $Y$ have a common  marginal $F$. 
However,  for other individuals, we just know $L=\min(X,Y)$  and we want to estimate $U=\max(X,Y)$ from $L$ (dependent censored data). To this purpose we can use  the regression curve $E(U|L=t)$ or the median regression curve obtained from the   conditional survival function $\Pr(U>y|L=x)$. 

To illustrate we simulate $m=100$ point with normal (Gaussian) and Exponential distributions for the common marginal $F$. 
Specifically, we choose $\mu=60$ and $\sigma=5$ in the normal model and  $\sigma=\mu=60$ in the exponential model. Then we assume different copula models. 

\begin{figure}[ptb]
	\begin{center}
		\includegraphics[scale=0.45]{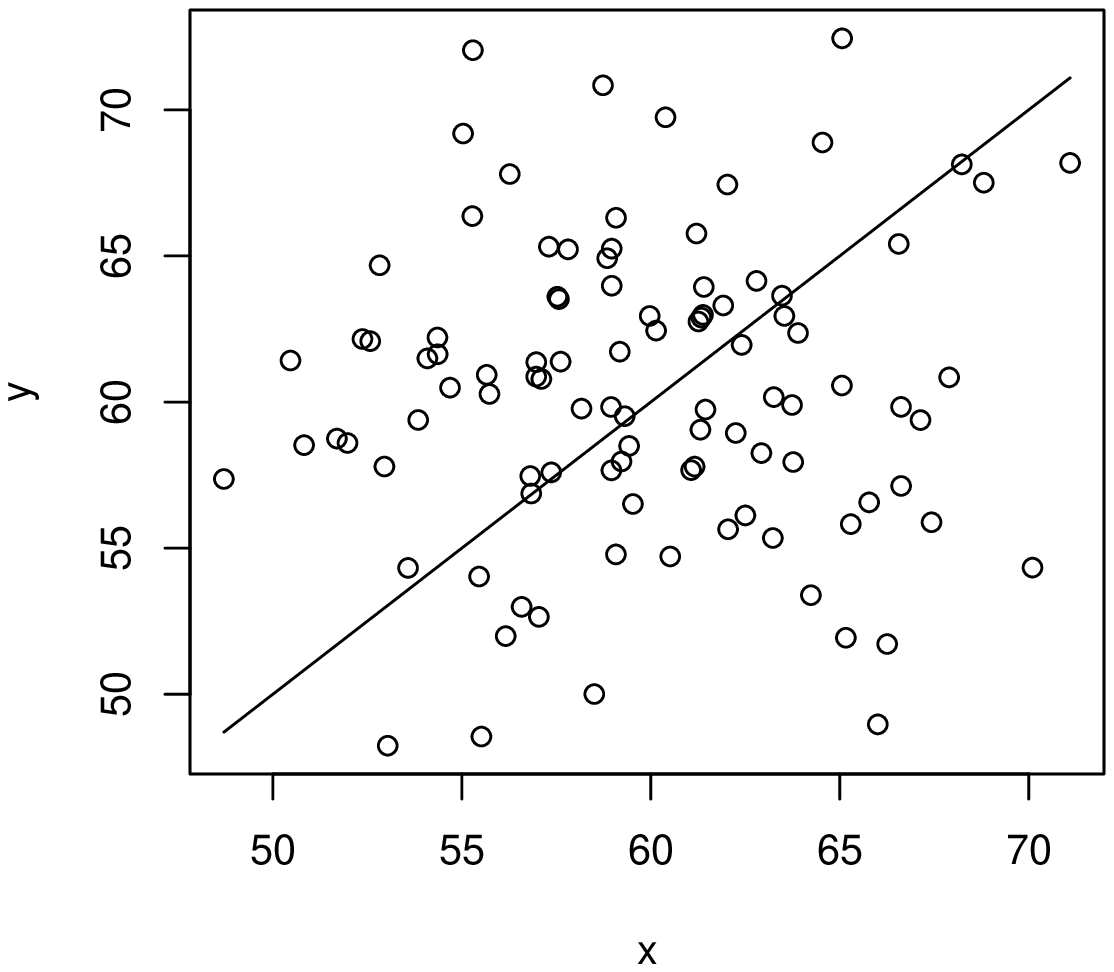}
		\includegraphics[scale=0.45]{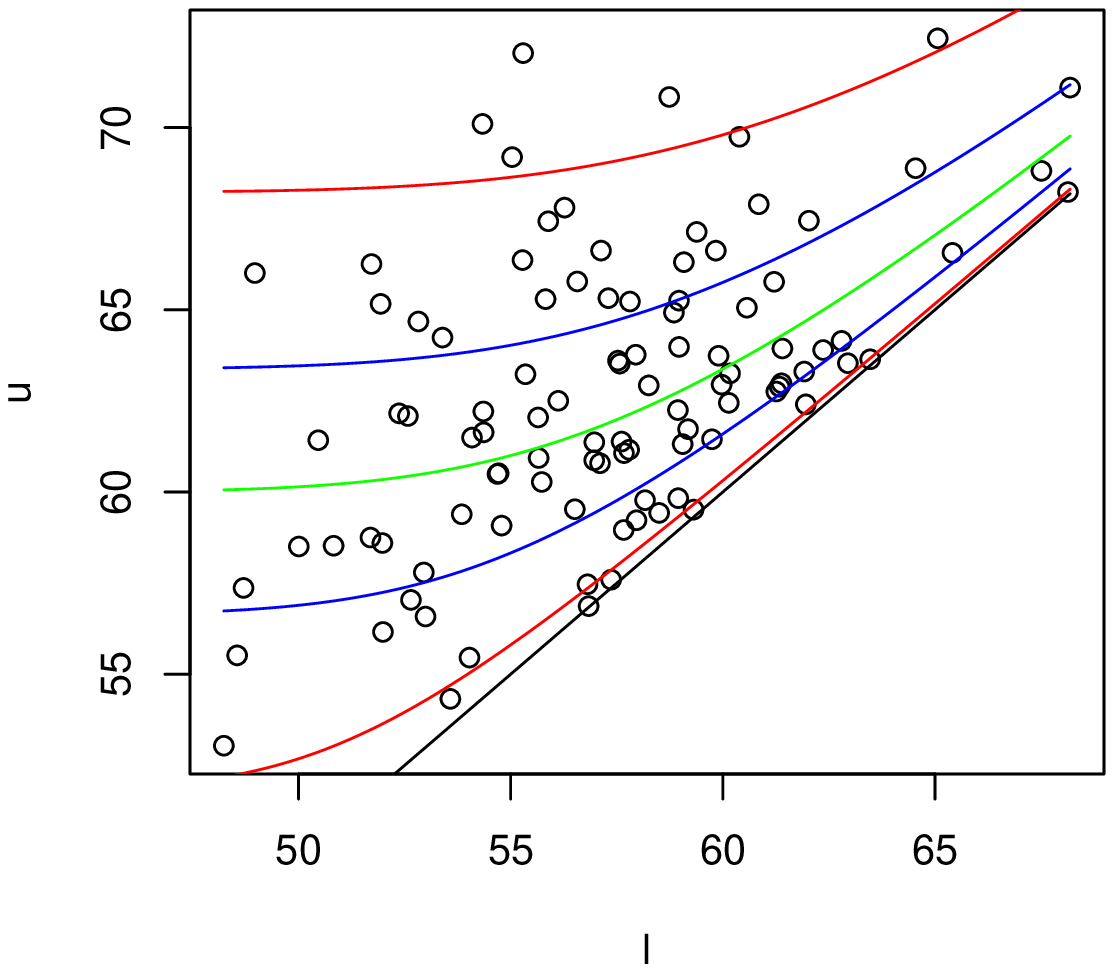}
		%470x470despacho
	\end{center}
	\caption{Plots of the data from $(X,Y)$ (left) and $(L,U)$ (right) when $X$ and $Y$ have independent normal distributions with  $\mu=60$ and $\sigma=5$. In the right plot we include the conditional median curve (green) and the  $90\%$ (red) and $50\%$ (blue) quantile-confidence  bands to estimate $U$ from $L$.}%
	\label{fig1}%
\end{figure}

In the first case, we want to note that $L$ and $U$ are dependent even if $X$ and $Y$ are independent. Thus we plot the data obtained from two independent normal distributions $N(60,5)$ in Figure \ref{fig1}, left, and the associated ordered data in Figure \ref{fig1}, right. The black lines represent the line $y=x$.  In the right plot we also include the conditional median curve (green) and the  $90\%$ (red) and $50\%$ (blue) quantile-confidence  bands computed from \eqref{IIDmedian}.   In this sample we have the following sample means for the ordered data: $\bar L=57.27662$, $\bar U=62.82893$ and $\bar U-\bar L=5.552306$.

\begin{figure}[ptb]
	\begin{center}
		\includegraphics[scale=0.45]{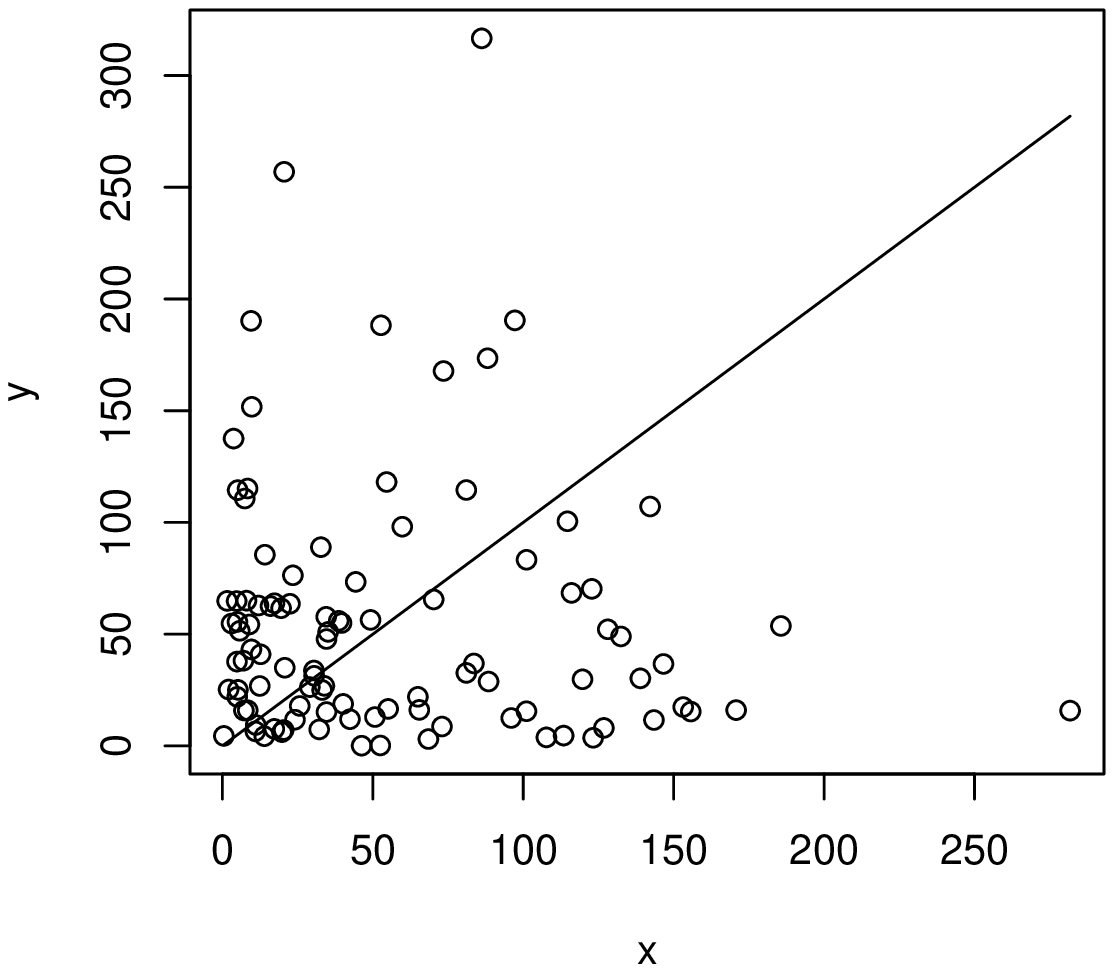}
		\includegraphics[scale=0.45]{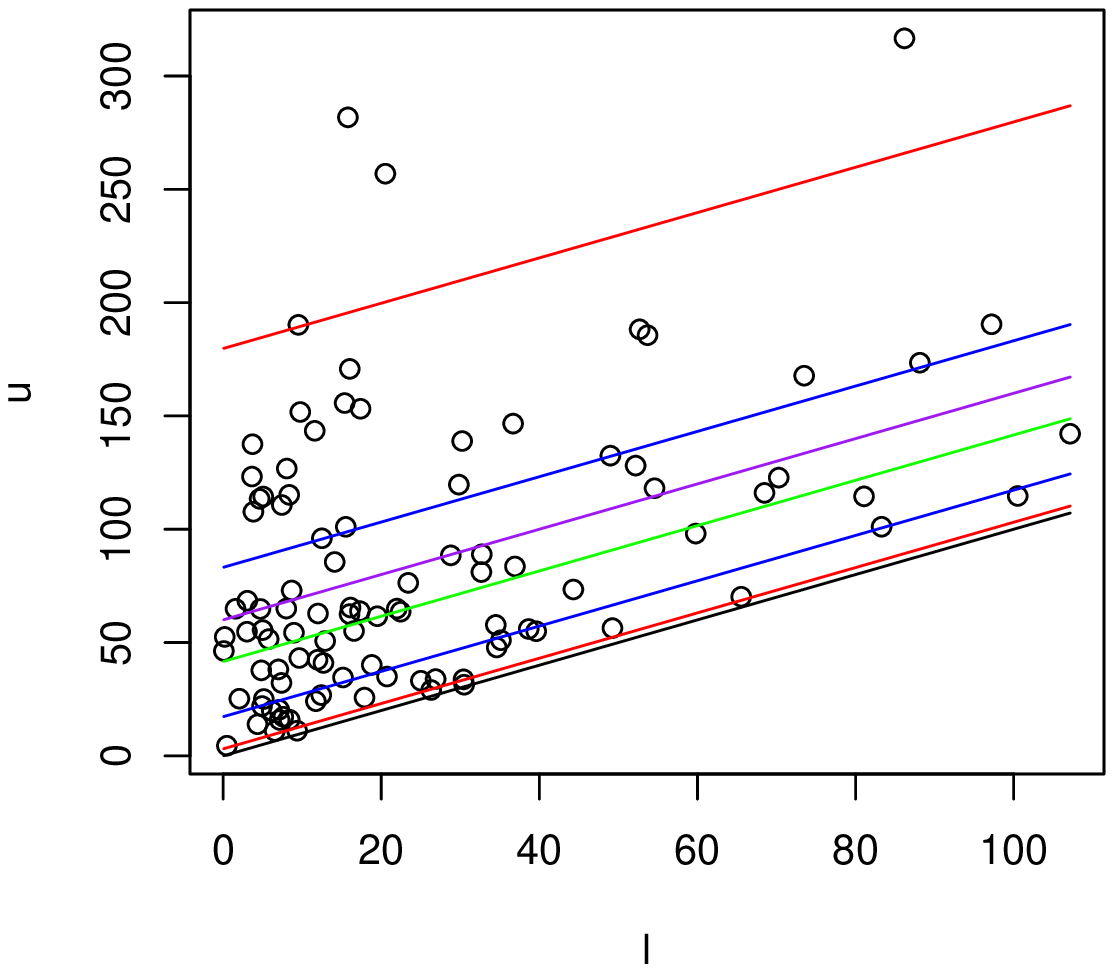}
		%470x470despacho
	\end{center}
	\caption{Plots of the data from $(X,Y)$ (left) and $(L,U)$ (right) when $X$ and $Y$ have independent exponential distributions with  $\mu=60$. In the right plot we include the regression line (purple), the conditional median curve (green) and the  $90\%$ (red) and $50\%$ (blue) quantile-confidence  bands to estimate $U$ from $L$.}%
	\label{fig2}%
\end{figure}

The analogous plots for the data obtained from two independent exponential distributions $Exp(60)$ are in  Figure \ref{fig2}. Note that  the conditional median curve (green) and the  $90\%$ (red) and $50\%$ (blue) quantile-confidence  bands  in the right plot are determined by lines as stated in \eqref{IIDmedianExp}. We also include the theoretical regression curve $E(U|L=x)=E(X|X>x)=60+x$ (purple line). Note that many of these data could be censored data in practice (e.g. all that are greater than 100) due to the big dispersion of the exponential model and the independence. In this sample we get $\bar L=25.35254$, $\bar U=84.22409$ and $\bar U-\bar V=58.87156$. Moreover, only two data from $L$ are greater than $100$ while we have  $35$ data from $U$ greater than $100$.  So the censure must be considered in practice.

Let us consider now that $(X,Y)$ are DID with a copula $C$ and a marginal distribution $F$. We consider again the above normal and exponential models for $F$. For $C$ we consider the following Clayton copula
\begin{equation}\label{C}
C(u,v)=\frac{uv}{u+v-uv}
\end{equation}
for $u,v\in[0,1]^2$ which represents a positive (symmetric) dependence between $X$ and $Y$. 

To generate a sample from $(X,Y)$, if $(U,V)$ is the random vector associated with this this copula,  we compute the conditional distribution  $(V|U=u)$ as 
$C_{2|1}(v|u)=\partial_1C(u,v)$
%=\frac{v^2}{(u+v-uv)^2}.$$
and its inverse function as 
$C^{-1}_{2|1}(q|u)=u/(u-1+q^{-1/2})$
for $u,q\in(0,1)$. Thus we can use the inverse transform method to get a sample from $C$. First, we generate a random (uniform) value $U$ in $(0,1)$ and then we use $C^{-1}_{2|1}(Q|U)$ to get the associated value $V$ from a random independent $Q\in(0,1)$. Then we use the transform $F^{-1}$ to obtain a sample from $(X,Y)$. This a standard method in copula theory (see, e.g., \cite{N06}, p.\ 41) . 

\begin{figure}[ptb]
	\begin{center}
		\includegraphics[scale=0.45]{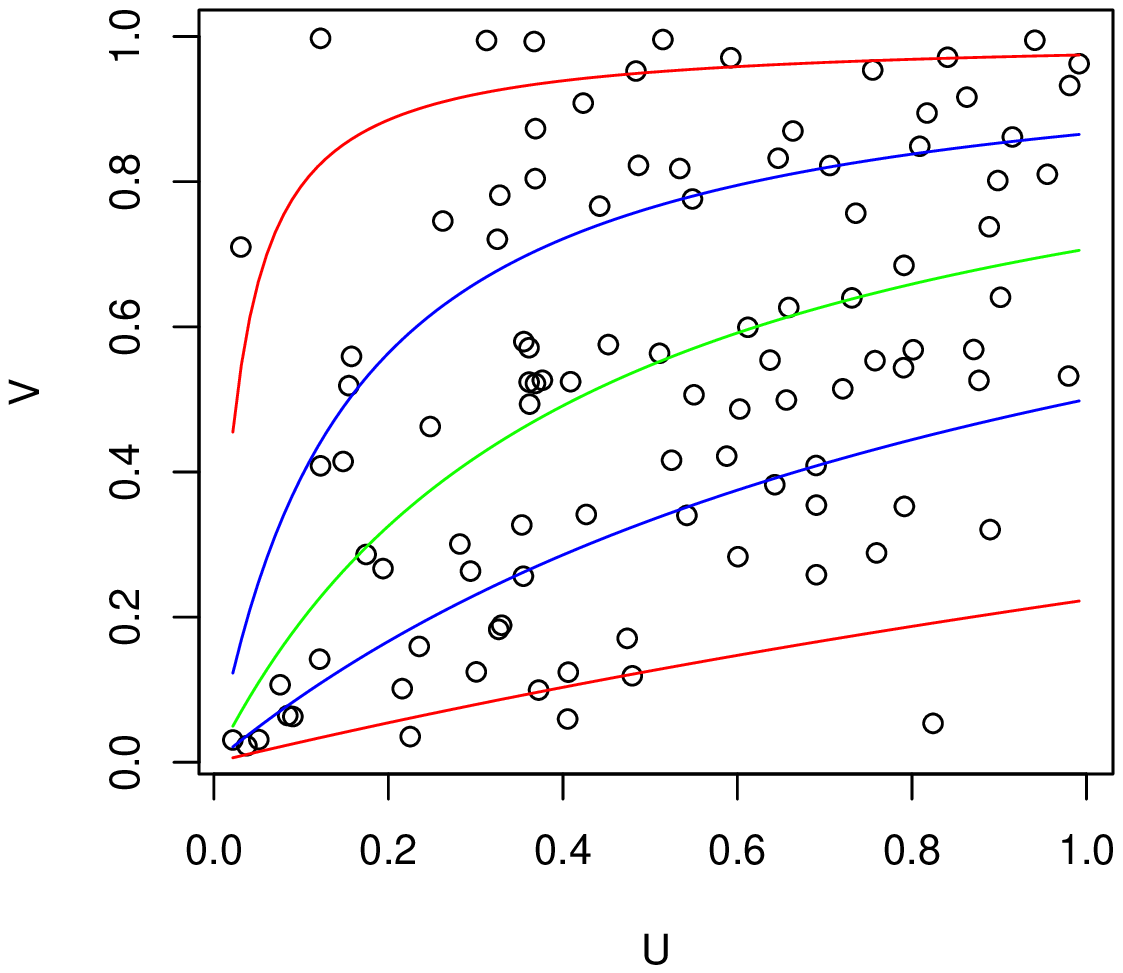}
		\includegraphics[scale=0.45]{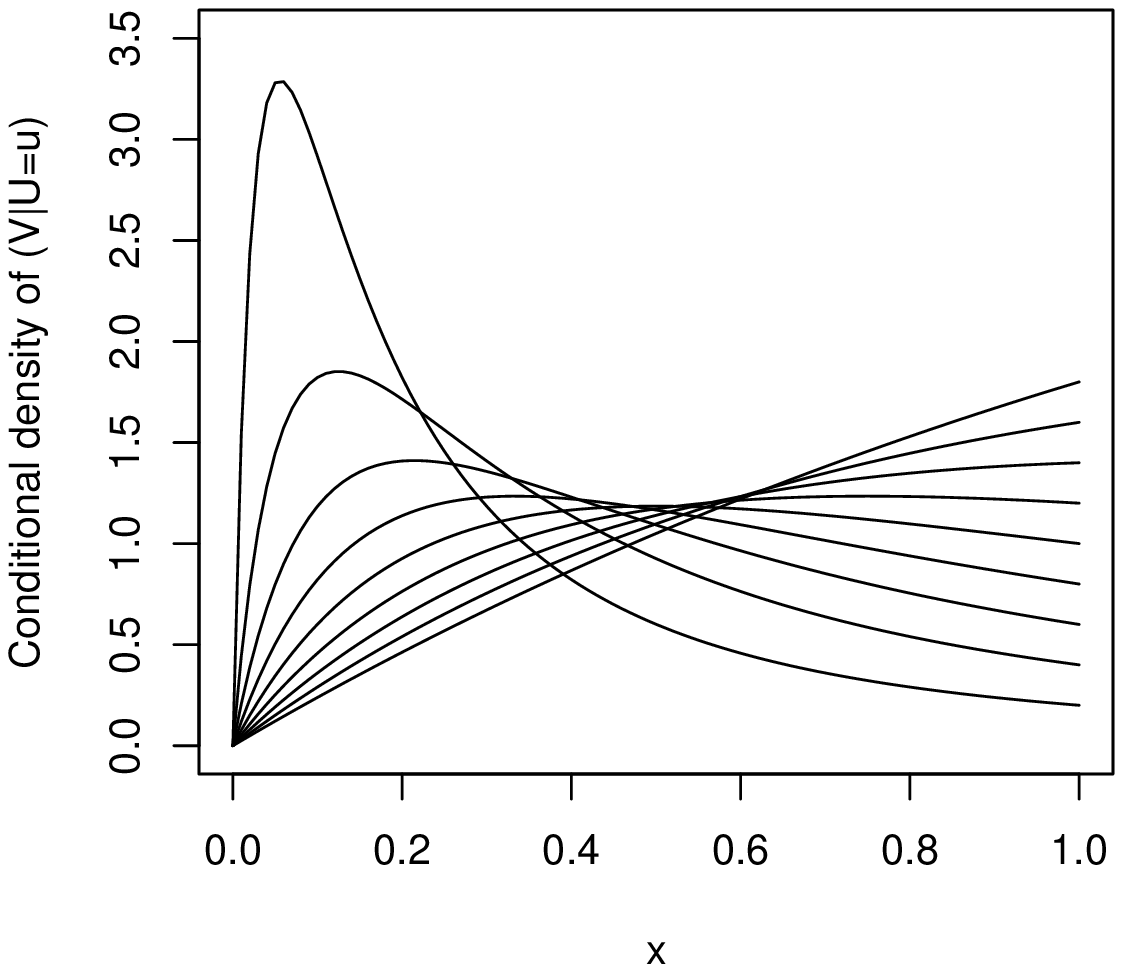}
		%470x470despacho
	\end{center}
	\caption{Plots  of the data from $(U,V)$ (left) from the Clayton copula $C$ in \eqref{C}  and conditional PDF (right) of $(V|U=u)$ for $u=0.1,\dots,0.9$ (from the top values at $x=0.1$).}%
	\label{fig3}%
\end{figure}

\begin{figure}[ptb]
	\begin{center}
		\includegraphics[scale=0.45]{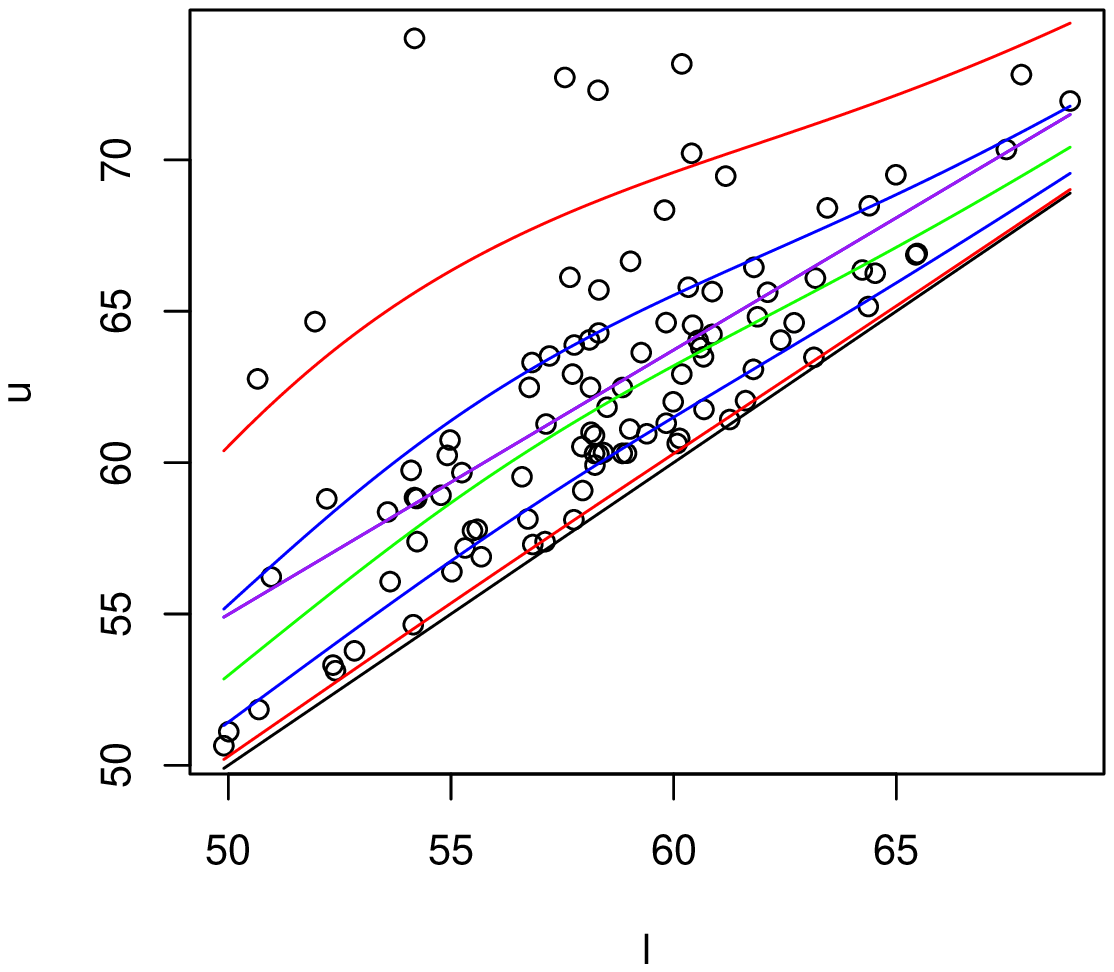}
		\includegraphics[scale=0.45]{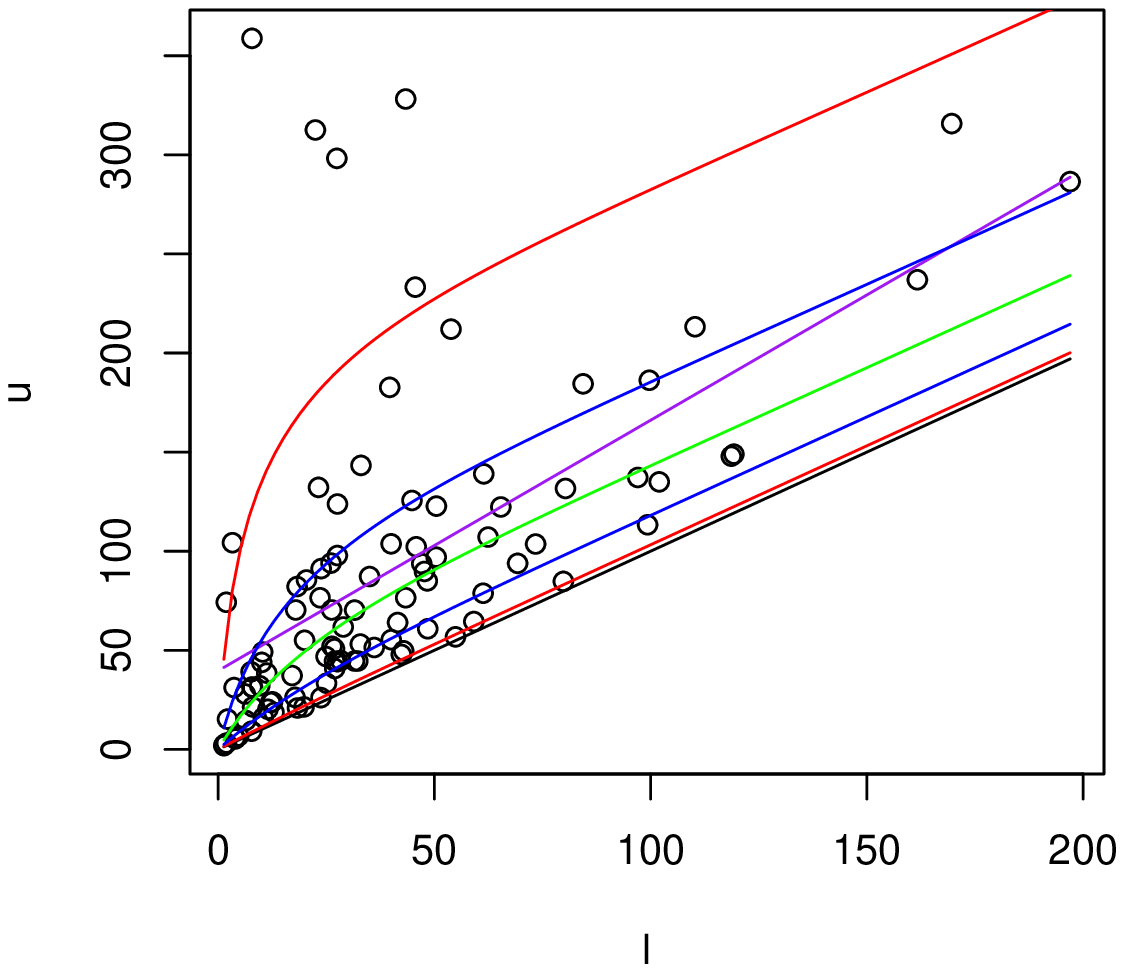}
		%470x470despacho
	\end{center}
	\caption{Plots of the data from  $(L,U)$  when $X$ and $Y$ have the Clayton copula $C$ in \eqref{C} and normal distributions with  $\mu=60$ and $\sigma=5$ (left) and exponential distributions with $\mu=60$ (right). We also include the regression line (purple), the conditional median curve (green) and the  $90\%$ (red) and $50\%$ (blue) quantile-confidence  bands.}%
	\label{fig4}%
\end{figure}

The data obtained from $C$ are plotted in Figure \ref{fig3}, left. Note that we can use  $C^{-1}_{2|1}$ to plot the conditional median curve and the $90\%$ and $50\%$ confidence bands for these data. In Figure \ref{fig3}, right, we plot the conditional PDF of $(V|U=u)$ for $u=0.1,\dots,0.9$. The data from  $(L,U)$ for the normal  and exponential  models are in Figure \ref{fig4}.

To get the conditional median curves and the  confidence bands for  $(L,U)$ we need the conditional distribution of $(U|L=x)$ that can be obtained 
from \eqref{G21} as
$
G_{2|1}(y|x)=D_{2|1} (F(y)|F(x))
$ 
for $y\geq x$, where 
$$D_{2|1}(v|F(x))=\frac{\partial_{1}D(F(x),v)}{\partial_{1}D(F(x),1)}$$
and
$$\partial_{1}D(u,v)=2\partial_{1}C(u,v)-2\partial_{1}C(u,u)=\frac{2v^2}{(u+v-uv)^2}-\frac{2}{(2-u)^2}$$
for $v\geq u$. Hence 
$$\partial_{1}D(u,1)=2-\frac{2}{(2-u)^2}=2\frac{u^2-4u+3}{(2-u)^2}.$$
To compute the inverse of $G_{2|1}$ we need to solve in $y$  the equation $G_{2|1}(y|x)=q$ for $q\in(0,1)$. This leads to 
%$$\partial_{1}D(F(x),F(y))=q \partial_{1}D(F(x),1)$$
%that is equivalent to 
$$\frac{F^2(y)}{(F(x)+F(y)-F(x)F(y))^2}=\frac{1-q+q(2-F(x))^2}{(2-F(x))^2}.$$
Therefore
%$$\frac{F(x)+F(y)-F(x)F(y)}{F(y)}=\sqrt{\frac{(2-F(x))^2}{1-q+q(2-F(x))^2}}$$
%and so we get
%$$\frac{F(x)}{F(y)}=F(x)-1+\frac{2-F(x)}{\sqrt{1-q+q(2-F(x))^2}}$$
%and
$$y=F^{-1}\left( \frac{F(x)}{F(x)-1+\frac{2-F(x)}{\sqrt{1-q+q(2-F(x))^2}}}\right).$$
We use this expression in Figure \ref{fig4}  to compute the conditional median curve and the associated $90\%$ and $50\%$ confidence bands for the above normal and exponential models. In these figures we also include the empirical regression lines to estimate $U$ from $L$ (purple lines). 
Note that in practice,  we will have a lot of censored data from the data in Figure \ref{fig4}, right. Actually, in this sample  we have $7$ values of $L$ and $31$ values of $U$ greater than $100$. If the data represent the ages (in years) for a disease in twins organs (e.g. breast cancer), the censure means that the patients died before they suffer this disease  in this organ (or that they do not have the disease when we finish the experiment).

%\begin{figure}[ptb]
%	\begin{center}
%	\includegraphics[scale=0.6]{Fig5l}	\includegraphics[scale=0.6]{Fig5r}
%470x470despacho
%\end{center}
%\caption{Plots of the data from $(X,Y)$ (left) and $(L,U)$ (right) when $X$ and $Y$ have exponential  distributions with  $\mu=60$ and the Clayton copula $C$ in \eqref{C}. We also include the regression line (purple), the conditional median curve (green) and the  $90\%$ (red) and $50\%$ (blue) quantile-confidence  bands.}%
%\label{fig5}%
%\end{figure}

\begin{figure}%[ptb]
	\begin{center}
\includegraphics[scale=0.45]{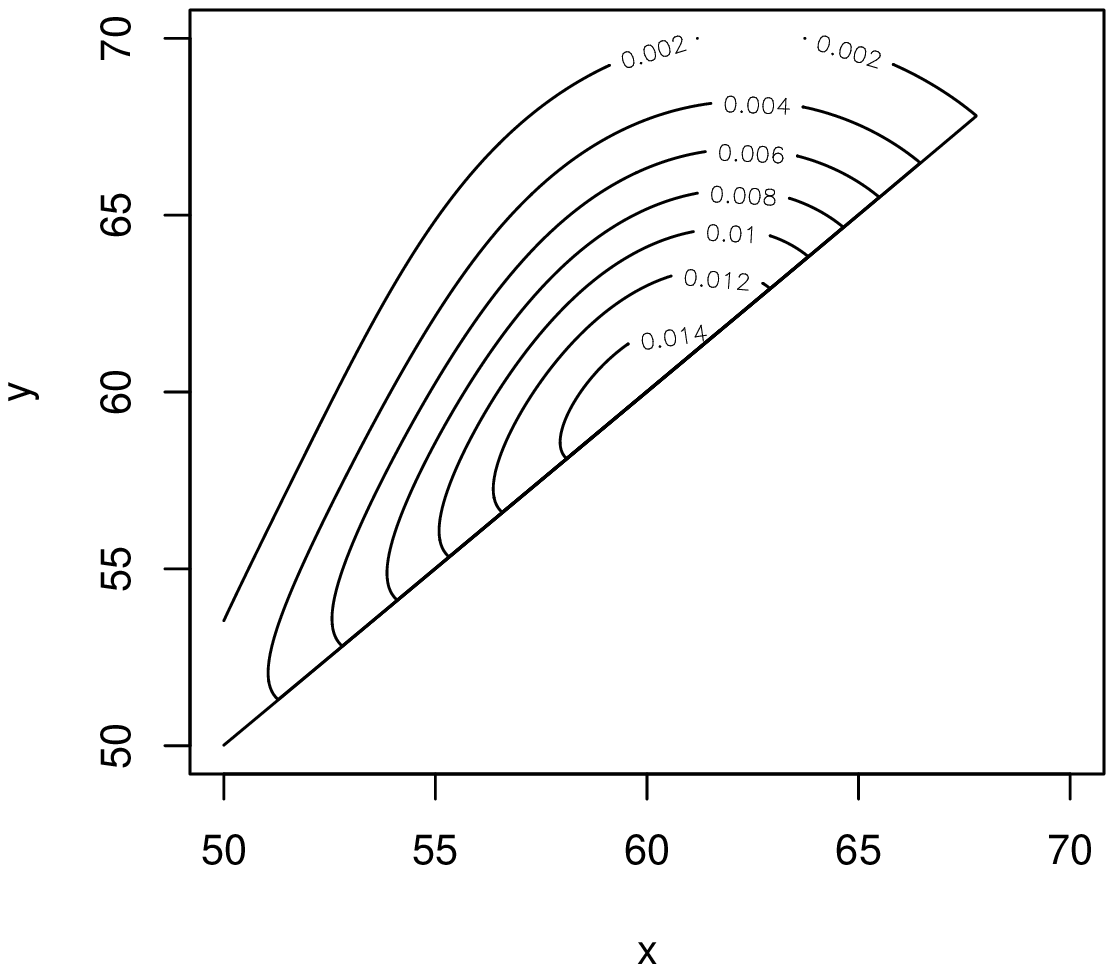}	\includegraphics[scale=0.45]{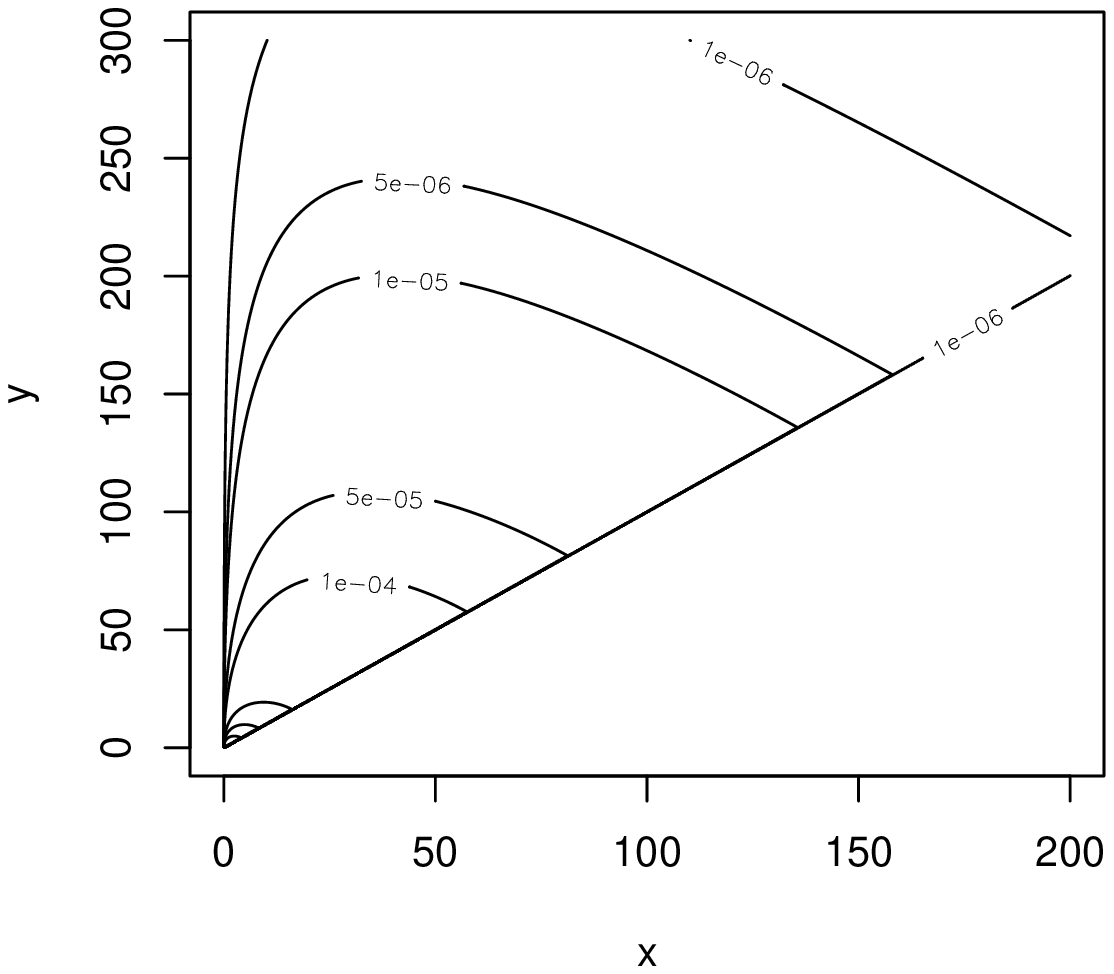}
		%470x470despacho
	\end{center}
	\caption{Contour plots for the joint PDF $\mathbf{g}$ of $(L,U)$  when $X$ and $Y$ have the Clayton copula $C$ in \eqref{C} and normal marginal distributions with $\mu=60$ and $\sigma=5$ (left) and exponential  distributions with  $\mu=60$ (right).}%
	\label{fig6}%
\end{figure}

\begin{figure}%[ptb]
	\begin{center}
\includegraphics[scale=0.45]{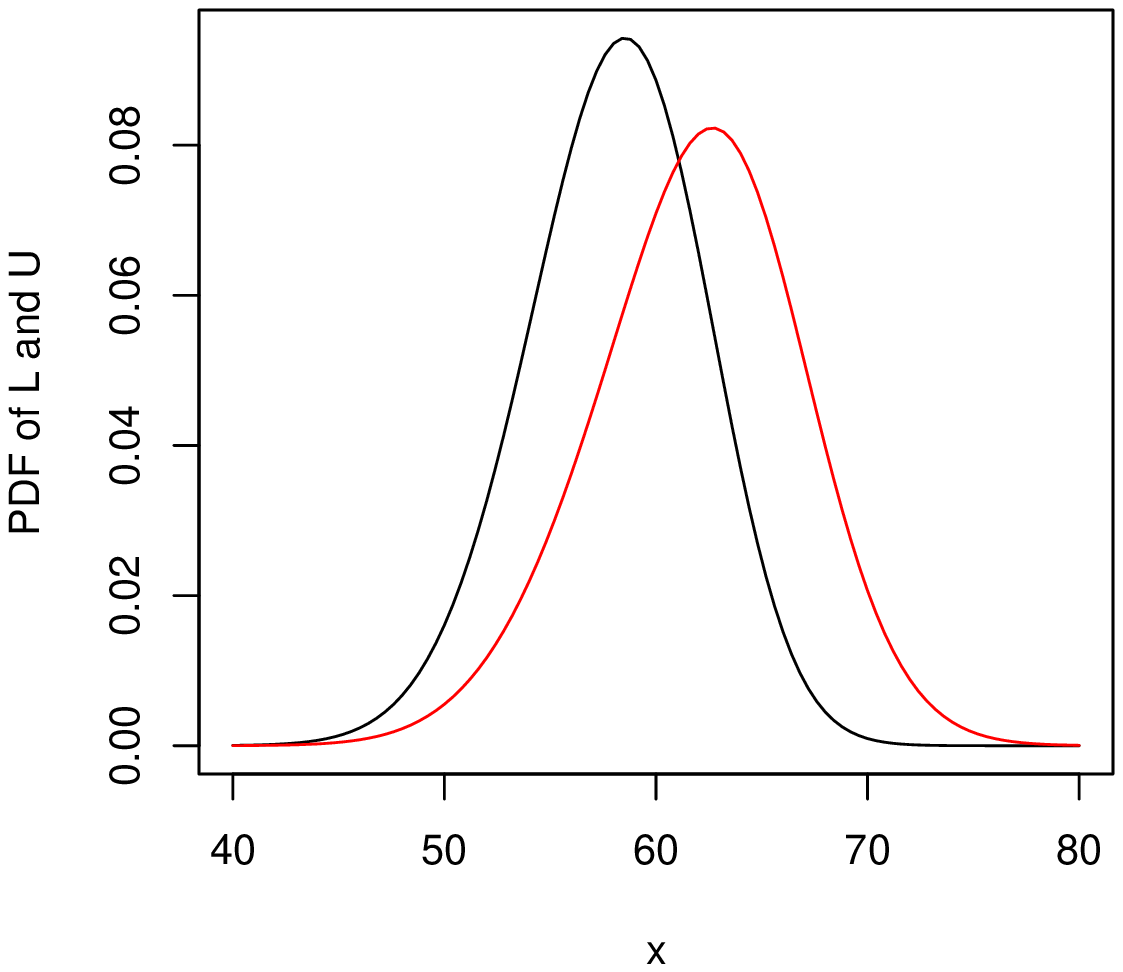}	\includegraphics[scale=0.45]{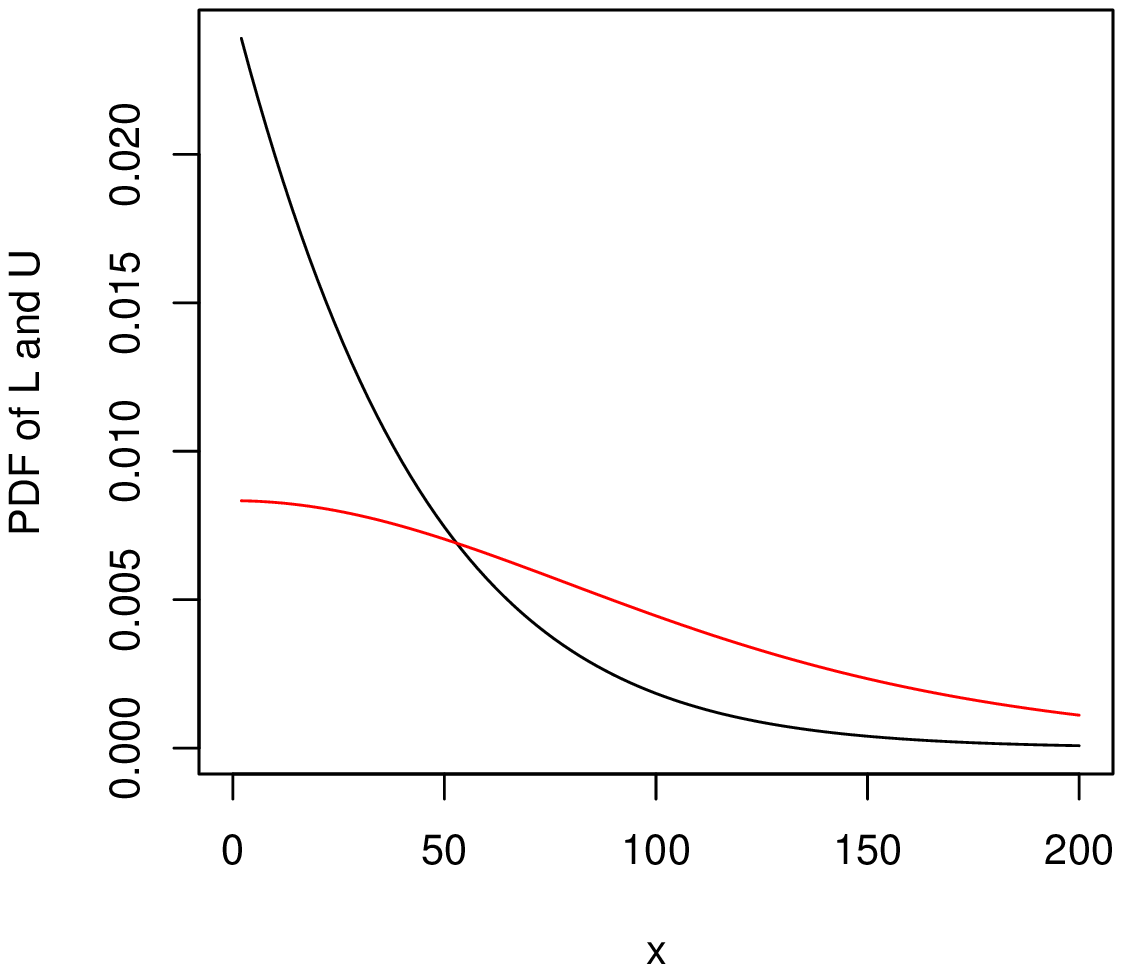}
		%470x470despacho
	\end{center}
	\caption{Plots for the PDF of $L$ (black) and $U$ (red) when $X$ and $Y$ have the Clayton copula $C$ in \eqref{C} and normal distributions with $\mu=60$ and $\sigma=5$ (left) and exponential  distributions with  $\mu=60$ (right).}%
	\label{fig7}%
\end{figure}

Finally, note that we can use the joint PDF $\mathbf{g}$ given in \eqref{g2} to plot the level curves for $(L,U)$. For the above Clayton (EXC)  copula we have
$$\mathbf{g}(x,y)=2f(x)f(y)c(F(x),F(y))$$
for $x\leq y$, where
$$c(u,v)=\partial_{1,2}C(u,v)=\frac{2uv}{(u+v-uv)^3}$$
for $u,v\in[0,1]$. For the above normal and exponential models we get the contour plots given in Figure \ref{fig6}. They fit to the data sets plotted in Figure \ref{fig4}. Analogously, we can plot the marginal PDF of $L$ and $U$ (see Figure \ref{fig7}). Note that for the above Clayton copula we obtain $D_1(u)=(3u-2u^2)/(2-u)$ and $D_2(u)=u/(2-u)$ for $u\in[0,1]$. So $L$ and $U$ do not have neither normal nor exponential distributions (but some plots are similar).

%\section{A real data example.}

\section{Conclusions}

The main purpose of this paper is to provide an alternative representation based on distortions   to the classical copula representation for the joint distribution of a random vector. The new representations are more flexible since they are not based on the marginals of the original model. However, a disadvantage is that they do not separate the dependence structure and the marginals.
The examples in Section 3 show that, in some cases, it could be better to use the new representation instead of the copula representation. 

The properties of the multivariate distorted distributions (MDD) are similar to that of copula representations. They allow us to get the marginal and conditional distributions. They can also be used to obtain the regression and median regression curves and the associated confidence bands.  Moreover, in some relevant examples studied in Section 3, we have unique MDD representations for fixed continuous distribution functions $G_1,\dots,G_n$. In these cases, the distortion function $D$ is uniquely determined by the copula $C$. 

We provide several examples where these representations are useful. More examples can be obtained in a similar way. We also include a simulation study for  paired ordered data from independent and dependent variables. This procedure can be used to predict the second failure in twins organs from the first one. It can also be used to complete (to estimate) censored data from dependent random variables.  

This paper is a first step and so it just contains the basic properties and some examples.  The main task for future research projects is to develop the appropriate inference procedures (and its properties) to apply these representations to real data sets. %This study must include the case of censored data.  
To this purpose we could use the wide literature about copula and regression quantiles estimations. Note that, in some relevant examples, the distortion function $D$ can be estimated from the copula $C$ and the training sample. The MDD representations can be also useful to study other relevant examples with dependent data.

%\section*{Acknowledgements}
%\acks

\begin{acknowledgements}
JN is partially supported by  Ministerio de Ciencia e Innovación  of Spain under grant PID2019-103971GB-I00, 
CC and ML  are partially supported by the GNAMPA research group of INDAM (Istituto Nazionale di Alta Matematica) and CC, ML and FD are also partially supported by MIUR-PRIN 2017, Project ``Stochastic Models for Complex Systems'' (No. 2017JFFHSH). 
\end{acknowledgements}

%\begin{acknowledgements}
%\end{acknowledgements}




\end{document}